\documentclass[a4paper,twoside,12pt]{article}
\usepackage{amsmath}
\usepackage{amsthm}
\usepackage{amsfonts}
\usepackage{amssymb, amscd}
\usepackage[dvips]{graphicx}
 \numberwithin{equation}{section}
\title{Quasi-hom-Lie Algebras, Central Extensions and 2-cocycle-like Identities}
\author{Daniel Larsson \\
    \footnotesize \emph{Centre for Mathematical Sciences, Department of
   Mathematics, Lund Institute}\\
   \footnotesize \emph{of Technology, Lund University, Box
   118, SE-221 00 Lund, Sweden}\\
   \footnotesize \emph{dlarsson@maths.lth.se,}\normalsize\and
   Sergei D. Silvestrov\\
   \footnotesize \emph{Centre for Mathematical Sciences, Department of
   Mathematics, Lund Institute}\\
   \footnotesize \emph{of Technology, Lund University, Box
   118, SE-221 00 Lund, Sweden}\\
   \footnotesize \emph{ sergei.silvestrov@math.lth.se}\normalsize}
\date{February 17, 2004}
\theoremstyle{definition}
\newtheorem{dfn}{Definition}
\newtheorem{example}{Example}
\newtheorem{remark}{Remark}

\theoremstyle{plain}
\newtheorem{prop}{Proposition}
\newtheorem{thm}[prop]{Theorem}

\newtheorem{lem}[prop]{Lemma}

\def\a{\alpha}
\def\b{\beta}

\def\s{\sigma }
\def\f{\varphi }
\def\e{\varepsilon }
\def\k{\Bbbk }
\def\o{\omega }

\def\i{\iota}
\def\C{\mathbb{C} }
\def\Z{\mathbb{Z} }

\def\mfD{\mathfrak{D} }
\def\g{\mathfrak{g} }
\def\gh{\hat{\mathfrak{g}} }
\def\gc{\check{\mathfrak{g}} }

\def\A{\mathcal A }
\def\U{\mathcal U }
\def\cbf{\mathbf{c}}
\def\De{ \Delta }
\def\d{ \delta }
\def\mfD{ \mathfrak{D} }
\def\M{ \mathfrak{M} }
\def\kbf{ \mathbf{k} }
\def\lbf{ \mathbf{l} }
\def\hbf{ \mathbf{h} }
\def\Witt{ \mathfrak{d} }

\def\cs{\circlearrowleft }
\def\afrak{\mathfrak{a} }

 \DeclareMathOperator{\id}{id}
 \DeclareMathOperator{\Mor}{Mor}

\DeclareMathOperator{\Ann}{Ann}\DeclareMathOperator{\Alt}{Alt}
\DeclareMathOperator{\pr}{pr}\DeclareMathOperator{\im}{im}\DeclareMathOperator{\sign}{sign}
\DeclareMathOperator{\Res}{Res}
\DeclareMathOperator{\Lin}{\mathcal{L}}
\begin{document}
\maketitle

\begin{abstract}
This paper begins by introducing the concept of a
\emph{quasi-hom-Lie algebra}, or simply, a \emph{qhl-algebra},
which is a natural generalization of hom-Lie algebras introduced
in a previous paper \cite{HartLarsSilv1C}. Quasi-hom-Lie algebras
include also as special cases (color) Lie algebras and
superalgebras, and can be seen as deformations of these by
homomorphisms, twisting the Jacobi identity and skew-symmetry. The
natural realm for these quasi-hom-Lie algebras is as a
generalization-deformation of the Witt algebra $\Witt$ of
derivations on the Laurent polynomials $\C[t,t^{-1}]$. We also
develop a theory of central extensions for qhl-algebras which can
be used to deform and generalize the Virasoro algebra by centrally
extending the deformed Witt type algebras constructed here. In
addition, we give a number of other interesting examples of
quasi-hom-Lie algebras, among them a deformation of the loop
algebra.
\end{abstract}
\footnotetext[1]{\emph{Keywords:} Color Lie algebras, Witt-type
algebras, deformation theory, $\sigma$-derivations, extensions,
Jacobi-type identities,
2-cocycle-like maps.\\
Mathematics Subject Classification 2000: 17B99 (Primary) 17B75,
17B68, 17A36, 17B40,
17B65,
17B66, 17B56 (Secondary)} \footnotetext[2]{The research was supported by the
Crafoord Foundation, the Swedish Royal Academy of Sciences and
Mittag-Leffler Institute.}
\newpage
\tableofcontents

\section{Introduction}
The classical Witt and Virasoro algebras are ubiquitous in
mathematics and theoretical physics, the latter algebra being the
unique one-dimens\-ional central extension of the former
\cite{QuantFieldStrings,DiFranMathiSenC,Fuchs1C,Fuchs2C,KacRainaC}.
Considering the origin of the Witt algebra this is not surprising:
the \emph{Witt algebra} $\Witt$ is the
  infinite-dimensional Lie algebra of complexified polynomial vector
  fields on the unit circle $S^1$. It can also be defined as
  $\Witt=\C\otimes\mathrm{Vect}(S^1)=\oplus_{n\in\Z}\C\cdot d_n$,
  where $d_n=-t^{n+1}d/dt$ is a linear basis for $\Witt$, and the Lie product being defined on the generators
  $d_n$ as
  $\langle d_n,d_m\rangle =(n-m)d_{n+m}$ and extended linearly to the whole $\Witt$.
This means in particular that any $\hat f\in\Witt$ can be written
as $\hat f=f\cdot d/dt$ with $f\in\C[t,t^{-1}]$, the algebra of
Laurent polynomials, and hence $\Witt$ can be viewed as the
(complex) Lie algebra of derivations on $\C[t,t^{-1}]$. When the
usual derivation operator is replaced by its difference
discretization or deformation, the underlying algebra is also in
general deformed, and the description and understanding of the
properties of the new algebra becomes a problem of key importance.

To put the present article into the right perspective and to see
where we are coming from we briefly recall the constructions from
\cite{HartLarsSilv1C}. In that paper we considered deformations of
$\Witt$ using $\sigma$-derivations, i.e. linear maps $D$
satisfying a generalized Leibniz rule $D(ab)=Da\cdot
b+\sigma(a)\cdot Db$. As we mentioned above the Witt algebra
$\Witt$ can be viewed as the Lie algebra of derivations on
$\C[t,t^{-1}]$. This observation was in fact our starting point in
\cite{HartLarsSilv1C} in constructing deformations of the Witt
algebra. Instead of just considering ordinary derivations on
$\C[t,t^{-1}]$ we considered $\sigma$-derivations. In fact, we did
something even more general as we considered a commutative,
associative $\C$-algebra $\A$ with $1_\A$ and a
$\sigma$-derivation $\Delta$ on $\A$. Forming the cyclic left
$\A$-module $\A\cdot\Delta$, a left submodule of the $\A$-module
$\mfD_\s(\A)$ of all $\sigma$-derivations on $\A$, we equipped
$\A\cdot\Delta$ with a bracket multiplication
$\langle\cdot,\cdot\rangle_\sigma$ such that it satisfied
skew-symmetry and a generalized Jacobi identity with six terms
\begin{align}\label{eq:jacobi_intro}
\cs_{x,y,z}\,\Big (\big\langle \sigma(x),\langle y,z\rangle_\s\big
\rangle_\s+\delta\cdot\big\langle x,\langle
y,z\rangle_\s\big\rangle_\s\Big )=0,
\end{align} where $\cs_{x,y,z}$ denotes cyclic summation
with respect to $x,y,z$ and where $\delta\in\A$. In the case when
$\A$ is a unique factorization domain we showed that the whole
$\A$-module $\mfD_\s(\A)$ is cyclic and can thus be generated by a
single element $\Delta$. Since $\C[t,t^{-1}]$ is a UFD this result
applies in particular to the $\sigma$-derivations on the Laurent
polynomials $\C[t,t^{-1}]$, and so we may regard
$\mfD_\s(\C[t,t^{-1}])$ as a deformation of
$\Witt=\mfD_{\id}(\C[t,t^{-1}])$. This means in particular that we
have a Jacobi-like identity (\ref{eq:jacobi_intro}) on
$\mfD_\s(\C[t,t^{-1}])$.

Furthermore, in \cite{HartLarsSilv1C} we concentrated mainly on
the case when $\delta\in\C\setminus\{0\}$ and so the Jacobi-like
identity (\ref{eq:jacobi_intro}) simplified to the Jacobi-like
identity with three terms
$$\cs_{x,y,z}\,\big\langle (\varsigma+\id)(x),\langle y,z\rangle_{\varsigma}\big
\rangle_{\varsigma}=0,$$ where
$\varsigma=\frac{1}{\delta}\bar\sigma$ is $1/\delta$-scaled
version of $\bar\sigma:\A\cdot\Delta\to\A\cdot\Delta$, acting on
this left module as
$\bar\sigma(a\cdot\Delta)=\sigma(a)\cdot\Delta$. Motivated by this
we called algebras with a three-term deformed Jacobi identity of
this form \emph{hom-Lie algebras}. Using that any algebra
\mbox{$\C$-endomorphism} $\s$ on $\C[t,t^{-1}]$ must be on the
form $\s(t)=qt^s$ for $s\in\Z$ and $q\in\C\setminus\{0\}$ unless
$\sigma(1)=0$ in which case $\sigma=0$ identically, we obtained a
$\Z$-parametric family of deformations which, when $s=1$, reduces
to a \mbox{$q$-deformation} of the Witt algebra becoming $\Witt$
when $q=1$. This deformation is closely related to the
$q$-deformations of the Witt algebra introduced and studied in
\mbox{\cite{AizSato1C,ChaiIsLukPopPresC,ChaiKuLukC,ChaiPres3C,ChungWSC,CurtrZachos1C,Kassel1C,PolychronakosC,SatoHT1C,SatoHT2C}.}
However, our defining commutation relations in this case look
somewhat different, as we obtained them, not from some conditions
aiming to resolve specifically the case of $q$-deformations, but
rather by choosing $\mathbb{C} [t,t^{-1}]$ as an example of the
underlying coefficient algebra and specifying $\sigma$ to be the
automorphism $\sigma_q: f(t)\mapsto f(qt)$ in our general
construction for $\sigma$-derivations. By simply choosing a
different coefficient algebra or basic $\sigma$-derivation one can
construct many other analogues and deformations of the Witt
algebra. The important feature of our approach is that, as in the
non-deformed case, the deformations and analogues of Witt algebra
obtained by various choices of the underlying coefficient algebra,
of the endomorphism $\s$ and of the basic $\s$-derivation, are
precisely the natural algebraic structures for the differential
and integral type calculi and geometry based on the corresponding
classes of generalized derivation and difference type operators.

We remarked in the beginning that the Witt algebra $\Witt$ has a
unique (up to isomorphism) one-dimensional central extension, the
\emph{Virasoro algebra}. In \cite{HartLarsSilv1C} we developed,
for the class of hom-Lie algebras, a theory of central extensions,
providing cohomological type conditions useful for showing the
existence of central extensions and for their construction. For
natural reasons we required that the central extension of a
hom-Lie algebra is also a hom-Lie algebra. In particular, the
standard theory of central extensions of Lie algebras becomes a
natural special case of the theory for hom-Lie algebras when no
non-identity twisting is present. In particular, this implies that
in the specific examples of deformation families of Witt and
Virasoro type algebras constructed within the framework of
\cite{HartLarsSilv1C}, the corresponding non-deformed Witt and
Virasoro type Lie algebras are included as the algebras
corresponding to those specific values of deformation parameters
which remove the non-trivial twisting. We rounded up
\cite{HartLarsSilv1C}, putting the central extension theory to the
test applying it for the construction of a hom-Lie algebra central
extension of the $q$-deformed Witt algebra producing a
$q$-deformation of Virasoro Lie algebra. For $q=1$ one indeed
recovers the usual Virasoro Lie algebra as is expected from our
general approach.

A number of examples of deformed algebras constructed in
\cite{HartLarsSilv1C} do not satisfy the three-term Jacobi-like
identity of hom-Lie algebras, but obey instead twisted six-term
Jacobi-like identities of the form (\ref{eq:jacobi_intro}). These
examples are recalled for the convenience of the reader among
other examples in Section 3. Moreover, there exists also many
examples where skew-symmetry is twisted as well. Taking the Jacobi
identity (\ref{eq:jacobi_intro}) as a stepping-stone we introduce
in this paper a further generalization of hom-Lie algebras by
twisting, not only the Jacobi identity, but also the skew-symmetry
and the homomorphism $\sigma$ itself ($\sigma$ is denoted by
$\alpha$ in this paper). In addition, we let go of the assumption
that $\delta$ (here denoted $\beta$) is an element of $\A$ and
assume instead that it is an endomorphism on $\A$. We call these
algebras \emph{quasi-hom-Lie algebras} or in short just
\emph{qhl-algebras}, see Definition 1. In this way we obtain a
class of algebras which not only includes hom-Lie algebras but
also color Lie algebras and other, more exotic types of algebras,
all of which can be viewed as a type of deformation of Lie
algebras. Note, however, that with these types of deformations we,
by design, leave the category of Lie algebras for a larger
category, including the Lie case as a subcategory.

The present paper is organized into two clearly distinguishable
parts. The first, consisting of Sections 2 and 3 concerns the
definition of qhl-algebras and some more or less elaborated
examples of such. The second part, \mbox{Section 4}, is devoted to
the (central) extension theory of qhl-algebras. Let us first
comment some on the first part. In Section 3 we give, based on
observations and results from \cite{HartLarsSilv1C}, examples of
qhl-algebras generalizations-deformations or analogues of the
classical Witt algebra $\Witt$, in addition to showing how the
notion of a qhl-algebra also encompasses Lie superalgebras and,
more generally, color Lie algebras by introducing gradings on the
underlying linear space and by suitable choices of deformation
maps. We also remark that we can define generalized color Lie
algebras by admitting the twists $\alpha$ and $\beta$. As another,
new, example of qhl-algebras we offer in subsection 3.2 a deformed
loop algebra. Section 4 is devoted to the development of a central
extension theory for qhl-algebras generalizing the theory for
(color) Lie algebras and hom-Lie algebras as developed in
\cite{HartLarsSilv1C}. We give necessary and sufficient conditions
for having a central extension and compare these results to the
ones given in the existing literature, for example
\cite{HartLarsSilv1C} for hom-Lie algebras and
\cite{Scheunert1C,ScheunertZhangC} for color Lie algebras. As a
last example we consider central extensions of deformed loop
qhl-algebras in subsection 4.3.

\section{Definitions and notations}
Throughout this paper we let $\k$ be a field of characteristic
zero and let $\Lin_\k(L)$ denote the linear space of $\k$-linear
maps (endomorphisms) of the $\k$-linear space $L$.
\begin{dfn}
A \emph{quasi-hom-Lie algebra} (\emph{qhl-algebra}, for short) is
a tuple \mbox{$(L,\langle \cdot,\cdot\rangle_L,\alpha,\beta,\o)$}
where
\begin{itemize}
    \item $L$ is a $\k$-linear space,
    \item $\langle\cdot,\cdot\rangle_L:L\times L\to L$ is a bilinear
      map
called a \emph{product} or
    a \emph{bracket in $L$}
    \item $\alpha,\beta:L\to L$, are linear maps,
    \item $\o:D_\o\to \Lin_\k(L)$ is a map with domain of definition $D_\o\subseteq L\times
    L$,
\end{itemize}
such that the following conditions hold:
\begin{itemize}
    \item ($\b$-twisting.) The map $\a$ is a $\b$-twisted
      algebra homomorphism,
      i.e. $$\langle\alpha(x),\alpha(y)\rangle_L=\beta\circ\alpha\langle x,y\rangle_L,
\quad\text{ for all } x,y\in L;$$
    \item ($\o$-symmetry.) The product satisfies a generalized skew-symmetry condition
        $$\langle x,y\rangle_L=\o(x,y)\langle y,x\rangle_L,
        \quad\text{ for all } (x,y)\in D_\o ;$$
    \item (qhl-Jacobi identity.) The bracket satisfies a generalized Jacobi identity
    $$\circlearrowleft_{x,y,z}\Big\{\,\o(z,x)\Big(\langle\alpha(x),\langle y,z\rangle_L\rangle_L+
    \beta\langle x,\langle y,z\rangle_L\rangle_L\Big)\Big\}=0,$$
     for all $(z,x),(x,y),(y,z)\in D_\o$.
\end{itemize}
\end{dfn}
\begin{remark}
    Note that if $\a=\id_L$ then, necessarily, $\b=\id|_{\langle L,L\rangle}$ restricted to the "commutator
    ideal"
     $\langle L,L\rangle\subseteq L$.
\end{remark}
\begin{remark}To avoid writing all the maps
  $\langle\cdot,\cdot\rangle_L,\,\a_L,\b_L$ and $\o_L$
every time we use two abbreviations for a qhl-algebra
$(L,\langle\cdot,\cdot\rangle_L, \a_L,\b_L,\o_L)$:
\begin{itemize}
    \item $(L,\Mor (L))$, "$\Mor$" for "morphism", or simply
    \item $L$, remembering that there are also maps present, implicitly, in the notation.
\end{itemize}
\end{remark}


 By a \emph{strong quasi-hom-Lie algebra morphism} (or simply a \emph{quasi-hom-Lie algebra morphism})
$$\phi:(L,\langle \cdot,\cdot\rangle_L,\alpha,\beta,\o)\to(L',\langle\cdot,\cdot\rangle_{L'},\alpha',\beta',\o')$$ we mean a linear map from $L$ to
$L'$ such that the following conditions hold for $x,y\in L$:
\begin{description}
    \item[M1.] $\phi(\langle x,y\rangle_L)=\langle\phi(x),\phi(y)\rangle_{L'}$,
    \item[M2.] $\phi\circ\alpha=\alpha'\circ\phi$,
    \item[M3.] $\phi\circ\beta=\beta'\circ\phi$.
\end{description}Conditions M2 and M3 are often referred to as "intertwining conditions".
\begin{lem}
    Let $\phi:L\to L'$ be a map satisfying condition \emph{M1}. Then
    $$\o_{L'}(\phi(x),\phi(y))\circ\phi=\phi\circ \o_L(x,y)$$ on
    $\langle L,L\rangle_L$ if $(x,y)\in D_{\o_L}$ and $(\phi(x),\phi(y))\in D_{\o_{L'}}$.
\end{lem}
\begin{proof}On the one hand
\begin{align*}
    \phi\langle x,y\rangle_L&=\langle\phi(x),\phi(y)\rangle_L=\o(\phi(x),\phi(y))\langle
    \phi(y),\phi(x)\rangle_L=\\
    &=\o(\phi(x),\phi(y))\circ\phi\langle y,x\rangle_L ,
\end{align*} and on the other,
\begin{align*}
    \phi\langle x,y\rangle_L=\phi\circ\o(x,y)\langle y,x\rangle_L.
\end{align*}Comparison proves the lemma.
\end{proof}

A \emph{weak quasi-hom-Lie algebra morphism} is a linear map $L\to
L'$ such that just condition M1 holds. If the weak (strong)
morphism $\f$ is a mono- and epimorphism it is a \emph{weak
(strong) isomorphism}.

\begin{dfn}
    Two qhl-algebras $L$ and $L'$ are said to be of \emph{$\f$-related $\o$-type} if there is a weak morphism $\f:L\to L'$ such
    that $$\o_{L'}(\f(x),\f(y))\circ\f=\f\circ\o_L(x,y)$$ for all $x,y\in L$.
\end{dfn}
\begin{dfn}
By an \emph{exact sequence} of qhl-algebras
$\{(L_i,\Mor(L_i))\}_{i\in\Z}$, with all $L_{i-1}$ and $L_i$ of
the same $\f$-related $\o$-type for each $i$, we mean a
commutative diagram with strong morphisms $\f_i$
\begin{align}\label{diag:exact_seq}
\begin{CD}
\cdots @>>> L_{i-1} @>\f_{i}>> L_i @>\f_{i+1}>> L_{i+1} @>>> \cdots \\
@. @V\a_{L_{i-1}} VV @V\a_{L_i}VV @V\a_{L_{i+1}} VV \\
\cdots @>>> L_{i-1} @>\f_{i}>> L_i @>\f_{i+1}>> L_{i+1} @>>> \cdots\\
@. @A\b_{L_{i-1}} AA @A\b_{L_i} AA @A\b_{L_{i+1}} AA\\
\cdots @>>> L_{i-1} @>\f_{i}>> L_i @>\f_{i+1}>> L_{i+1} @>>> \cdots\\
\end{CD}
\end{align}such that the rows are exact, that is, $\im (\f_i)=\ker (\f_{i+1})$.
\end{dfn}
\begin{dfn}
A \emph{short exact sequence} of qhl-algebras $L_i$ is a
commutative diagram as {\rm (\ref{diag:exact_seq})} with all
algebras but three consecutive ones zero, i.e. a diagram
\begin{align}\label{diag:ses_ab1}
\begin{CD}
0 @>>> \afrak @>\iota>> E @>\pr>> L @>>> 0\\
@. @V\a_{\afrak} VV @V\a_E VV @V\a_L VV \\
0 @>>> \afrak @>\iota>> E @>\pr>> L @>>> 0\\
@. @A\b_{\afrak} AA @A\b_E AA @A\b_L AA \\
0 @>>> \afrak @>\iota>> E @>\pr>> L @>>> 0\\
\end{CD}
\end{align}with commutative ''boxes'' and exact rows where $\i$ and $\pr$ are strong morphisms.
\end{dfn}
\begin{remark} Note that this means, in particular, that
    \begin{itemize}
        \item $\iota \circ
        \o_\afrak(a,b)=\o_E(\iota(a),\iota(b))\circ\iota$.
        \item $\pr\circ\o_E(e,e')=\o_L(\pr(e),\pr(e'))\circ\pr$.
    \end{itemize}
\end{remark}
\begin{dfn}
    A short exact sequence as {\rm(\ref{diag:ses_ab1})} is a \emph{quasi-hom-Lie algebra extension}
    of $L$ by $\afrak$, or by a slight abuse of language, we say that
    $E$
is an \emph{extension} of $L$ by $\afrak$.
\end{dfn}
\section{Examples}
\begin{example}\label{ex:homLie}
    By taking $\b$ to be the identity $\id_L$ and $\o=-\id_L$ we get
    the hom-Lie algebras discussed in a previous paper {\rm\cite{HartLarsSilv1C}}.
    We recall the definition for the reader's convenience. A
    \emph{hom-Lie algebra} is a non-associative algebra $L$ with
    bracket multiplication $\langle\cdot ,\cdot\rangle_\alpha$ with $\alpha:L\to
    L$ an algebra endomorphism, such that
    \begin{itemize}
        \item $\langle x,y\rangle_\alpha=-\langle
        y,x\rangle_\alpha,\qquad\qquad\qquad\qquad\,\,$ (Skew-symmetry)
        \item $\cs_{x,y,z}\, \big\langle (\id_L+\alpha)(x),\langle
        y,z\rangle_\alpha\big\rangle_\alpha=0\qquad$ ($\alpha$-deformed Jacobi
        identity).
    \end{itemize}
        Restricting further we get a Lie algebra by taking $\a$ also equal to the identity $\id_L$.\qed
\end{example}
Let $\Gamma$ be an abelian group. Then a \emph{$\Gamma$-graded
linear space} is a linear space $V$ with a direct sum
decomposition into subspaces labeled by $\Gamma$, i.e.
$$V=\bigoplus_{\gamma\in\Gamma}V_\gamma.$$ This means that any element in $V$ can be written uniquely as a \emph{finite} sum
$v=\sum_{\gamma\in \Gamma}v_\gamma$. The elements $v_\gamma\in
V_\gamma$ are called \emph{homogeneous of degree
$\gamma$}.\\
 A \emph{$\Gamma$-graded algebra} is a $\Gamma$-graded linear space with bilinear multiplication $\star$
respecting the grading in the following sense
$$V_{\gamma_1}\star V_{\gamma_2}\subseteq V_{\gamma_1+\gamma_2}.$$
\begin{example}\label{ex:superLie} A \emph{Lie superalgebra} is a $\Z_2$-graded ($\Z_2=\Z/2\Z$) algebra
$L$ decomposing as $L=L_0\oplus L_1$ with a bilinear product
$\langle\cdot,\cdot\rangle$ satisfying
\begin{itemize}
    \item $\langle L_{\gamma_1},L_{\gamma_2}\rangle\subseteq L_{\gamma_1+\gamma_2}$
    \item $\langle x,y\rangle=-(-1)^{\gamma_x\gamma_y}\langle y,x\rangle$ for $x\in L_{\gamma_x}$ and $y\in L_{\gamma_y}$, and
    \item $(-1)^{\gamma_z\gamma_x}\big \langle x,\langle y,z\rangle\big \rangle+(-1)^{\gamma_x\gamma_y}
    \big \langle y,\langle z,x\rangle \big \rangle+
    (-1)^{\gamma_y\gamma_z}\big \langle z,\langle x,y\rangle \big
    \rangle=0$, for $x\in L_{\gamma_x}$, $y\in L_{\gamma_y}$ and $z\in L_{\gamma_z}$.
    \end{itemize} Note that $L_0$ is a
    subalgebra of $L$ but that $L_1$ is not. $L_0$ is the \emph{even} component of $L$ and $L_1$ the \emph{odd}.
    Lie superalgebras are qhl-algebras. This can be easily seen letting $L$ be $\Z_2$-decomposed in the definition of the
    qhl-algebras as $L=L_0\oplus L_1$, $\alpha=\beta=\id_L$ and
    $\o(x,y)v=\o(\gamma_x,\gamma_y)v=-(-1)^{\gamma_x\gamma_y}v$
    for $v\in L$ and where $(x,y)\in D_\o=L_0\cup L_1$ and $\gamma_x,\gamma_y\in\Z_2$ are the graded
    degrees (or \emph{parity}) of $x$ and $y$ respectively. \qed
\end{example}
\begin{example}\label{ex:colorLie}
    The super Lie algebra example is covered by the more general notion of a \emph{color Lie algebra} or $\Gamma$-graded
    $\varepsilon$-Lie algebra. In this case $\Gamma$ is any abelian group and the color Lie algebra $L$ with bracket $\langle\cdot,\cdot\rangle$ decomposes as
$L=\oplus_{\gamma\in\Gamma} L_\gamma$ where $\langle
L_{\gamma_1},L_{\gamma_2}\rangle\subseteq L_{\gamma_1+\gamma_2}$,
for $\gamma_1,\gamma_2\in\Gamma$. In addition, the "color
structure" includes a map $\e:\Gamma\times\Gamma\to \k$, called a
\emph{commutation factor}, satisfying
\begin{itemize}
    \item $\e(\gamma_x,\gamma_y)\e(\gamma_y,\gamma_x)=1$,
    \item $\e(\gamma_x+\gamma_y,\gamma_z)=\e(\gamma_x,\gamma_z)\e(\gamma_y,\gamma_z),$
    \item $\e(\gamma_x,\gamma_y+\gamma_z)=\e(\gamma_x,\gamma_y)\e(\gamma_x,\gamma_z)$,
\end{itemize}for $\gamma_x,\gamma_y,\gamma_z\in\Gamma$. The color skew-symmetry and Jacobi condition are now
stated, with the aid of $\e$, as
\begin{itemize}
    \item $\langle x,y\rangle=-\e(\gamma_x,\gamma_y)\langle y,x\rangle$
    \item $\e(\gamma_z,\gamma_x)\big \langle x,\langle y,z\rangle \big \rangle+\e(\gamma_x,\gamma_y)\big \langle y,
    \langle z,x\rangle \big \rangle+
    \e(\gamma_y,\gamma_z)\big \langle z,\langle x,y\rangle \big \rangle=0$
\end{itemize}for $x\in L_{\gamma_x}, y\in L_{\gamma_y}$ and $z\in L_{\gamma_z}$.
Color Lie algebras are examples of \mbox{qhl-algebras}. This can
be seen by putting a grading on $L$ in the definition of
qhl-algebras $L=\oplus_{\gamma\in\Gamma} L_\gamma$, $\a=\b=\id_L$
and $\o(x,y)v=-\e(\gamma_x,\gamma_y)v$, for $v\in L$ and where
$(x,y)\in D_\o=\cup_{\gamma\in\Gamma}L_\gamma$ and
$\gamma_x,\gamma_y\in\Gamma$ are the graded degrees of $x$ and
$y$. The $\o$-symmetry and the qhl-Jacobi identity gives us the
respective identities in the definition of a color Lie algebra.
\qed
\end{example}
\begin{remark}Since $\alpha_L=\beta_L=\id_L$ for (color) Lie
algebras,
  there is only one notion of morphism in this case, namely the usual
  (color) Lie algebra homomorphism.
\end{remark}
\begin{remark}
       By not restricting $\alpha_L$ to be the identity in Example \ref{ex:colorLie} we can define
       \emph{color hom-Lie algebras}, and similarly, with
       $\beta_L\neq\id_L$, \emph{color qhl-algebras}.
\end{remark}

In the next section we summarize several results from
\cite{HartLarsSilv1C} providing non-trivial examples of
qhl-algebras. For proofs and more details we refer to
\cite{HartLarsSilv1C}.

\subsection{$\sigma$-derivations}
In this section, we let $\A$ denote a commutative, associative
$\k$-algebra with unity, and let $\mfD_\s(\A)$ denote the set of
$\sigma$-derivations on $\A$, that is the set of all $\k$-linear
maps $D:\A\to\A$ satisfying the $\s$-Leibniz rule
$$D(ab)=D(a)b+\s(a)D(b).$$

We now fix a homomorphism $\s:\A\to\A$, an element
$\De\in\mfD_\s(\A)$, and an element $\d\in\A$, and we assume that
these objects satisfy the following two conditions:
\begin{align}\label{eq:GenWittCondC1}
    \s(\Ann(\De))\subseteq \Ann(\De),
\end{align}
\begin{align}\label{eq:GenWittCondC2}
    \De(\s(a)) = \d\s(\De(a)),\quad\text{for }a\in\A,
\end{align} where $\Ann(\Delta)=\{a\in\A\,|\, a\cdot
\Delta=0\}$. Let
$$\A\cdot\De=\{a\cdot\De\;|\;a\in\A\}$$
denote the cyclic $\A$-submodule of $\mfD_\s(\A)$ generated by
$\De$ and extend $\s$ to $\A\cdot\De$ by
$\sigma(a\cdot\De)=\sigma(a)\cdot\De$. We have the following
theorem, which introduces a $\k$-algebra structure on
$\A\cdot\De$.
\begin{thm}\emph{(cf. \cite{HartLarsSilv1C})} \label{thm:GenWittC}
If {\rm (\ref{eq:GenWittCondC1})} holds then the map
\begin{align*}
    \langle\cdot,\cdot\rangle_\s:\A\cdot\De\times\A\cdot\De \to
    \A\cdot\De
\end{align*} defined by setting
\begin{align} \label{eq:GenWittProdDefC}
    \langle a\cdot\De,b\cdot\De\rangle_\s=(\s(a)\cdot\De)\circ(b\cdot\De)-(\s(b)\cdot\De)
    \circ(a\cdot\De), \quad\text{for }a,b\in\A,
\end{align}
where $\circ$ denotes elementwise composition, is a well-defined
$\k$-algebra product on the $\k$-linear space $\A\cdot\De$, and it
satisfies the following identities:
\begin{align}\label{eq:GenWittProdFormulaC}
    \langle a\cdot\De, b\cdot\De\rangle_\s=\big(\s(a)\De(b)-\s(b)\De(a)\big)\cdot\De,
\end{align}
\begin{align}\label{eq:GenWittSkewC}
    \langle a\cdot\De, b\cdot\De\rangle_\s=-\langle b\cdot\De,
    a\cdot\De\rangle_\s,
\end{align}for $a,b,c\in\A$. In addition, if {\rm (\ref{eq:GenWittCondC2})} holds, then
\begin{equation} \label{eq:GenWittJacobiC}
    \cs_{a,b,c}\,\Big(\big \langle\s(a)\cdot\De,\langle
    b\cdot\De,c\cdot\De
\rangle_\s\big\rangle_\s+\d\cdot\big\langle a\cdot\De, \langle
b\cdot\De,c\cdot\De\rangle_\s\big\rangle_\s\Big)=0.
\end{equation}
\end{thm}The algebra $\A\cdot\De$ from the theorem is then a
qhl-algebra with $\a=\s$, $\b=\d$ and $\o=-1$.
\begin{remark}
     Let $\Delta$ be a non-empty family of commuting $\sigma$-derivations on $\A$
     closed under composition of maps. Then $\Delta$ generates a
     left $\A$-module $\A\otimes\Delta$ via the rule $b(a\otimes d)=(ba)\otimes d$,
     where $a,b\in\A$ and $d\in\De$. We extend any $d\in\Delta$ from $\A$ to
     $$\A\otimes\Delta=\{\sum_{\text{finite}} a\otimes d\,|\,a\in\A,d\in\Delta\}$$
     by the rule
     $$d(a\otimes d')=d(a)\otimes d'+\sigma(a)\otimes dd',$$ where $dd'$
     denotes (associative) composition $dd'(a)=d(d'(a))$. For $a\in\A$ and $d\in\Delta$
     we can identify $a\otimes d$ and $a d$ as operators on $\A$ by $a\otimes d(r)=a(d(r))$
     for $r\in\A$.
     Define a product on monomials of $\A\otimes\Delta$ by
     $$\langle
     a\otimes d_1,b\otimes d_2\rangle_\sigma=\sigma(a)\otimes d_1(b\otimes d_2)-
     \sigma(b)\otimes d_2
     (a\otimes d_1),$$ and extend linearly to the whole $\A\otimes\Delta$. Then a
simple calculation using the commutativity of $\A$ and $\Delta$
shows that
     $$\langle
     a\otimes d_1,b\otimes d_2\rangle_\sigma=\big(\sigma(a) d_1(b)\big )\otimes d_2-
     \big(\sigma(b) d_2(a)\big )
     \otimes d_1.$$ Skew-symmetry also follows from this.
     Note also that if $d_1,d_2\in\Delta$ then \mbox{$d_1-d_2\in\mfD_\s(\A)$.} If $\Delta$ is maximal
     with respect to being commutative then
     \mbox{$d_1-d_2\in\Delta$.} We see that
     part of the above theorem generalizes to a setting with multiple $\sigma$-derivations. However, if there
     is a nice Jacobi-like identity as in the theorem is uncertain at this moment.
     The
     above construction parallels the one given in {\rm\cite{Passman1C}} with the difference that
     {\rm\cite{Passman1C}}
     considers the construction in a color Lie algebra setting.
\end{remark}

Under the assumption that $\A$ is a unique factorization domain
there exists $\De\in\mfD_\s(\A)$ such that $\A\cdot\De =
\mfD_\s(\A)$. This means that $\mfD_\s(\A)$ affords a qhl-algebra
structure under the above assumptions on $\A$. For more details
see \cite{HartLarsSilv1C}.

We now apply these ideas for the construction of some explicit
examples.

\subsubsection{A $q$-deformed Witt algebra}
Let $\A$ be the unique factorization domain $\k[t,t^{-1}]$, the
algebra of Laurent polynomials in $t$ over the field $\k$. Then
$\mfD_\sigma(\A)$ can be generated by a single element $D$ as a
left $\A$-module, that is $\mfD_\sigma(\A)=\A\cdot D$. When
$\sigma(t)=qt$ and $q\neq 1,0$, one can take
$$D=\frac{\id-\sigma}{1-q}\quad :\quad f(t)\mapsto\frac{f(t)-f(qt)}{1-q}=\frac{f(qt)-f(t)}{q-1}.$$
Then (\ref{eq:GenWittCondC1}) and (\ref{eq:GenWittCondC2}) are
satisfied with $\delta=1$, and Theorem \ref{thm:GenWittC} yields
the following qhl-algebra structure on $\mfD_\s(\A)$. This is, in
fact, an example of a qhl-algebra which is also a hom-Lie algebra
(see Example 1).
\begin{thm}\emph{(cf. \cite{HartLarsSilv1C})}\label{th:qWittlinearC} Let $\A=\k[t,t^{-1}]$. Then the $\k$-linear space
    $$\mfD_\s (\A) = \bigoplus_{n\in\Z} \k\cdot d_n,$$ with
    $d_n=-t^nD$ can be equipped with the
    bracket multiplication
    $$\langle\,\cdot,\,\cdot\rangle_\s\,:\,
\mfD_\s (\A)\times\mfD_\s (\A) \longrightarrow
    \mfD_\s (\A)$$ defined on generators by {\rm (\ref{eq:GenWittProdDefC})} as
    $$\langle d_n,d_m\rangle_\s=q^nd_{n}d_m-q^md_{m}d_n$$ with commutation relations
    $$\langle d_n,d_m\rangle_\s=\big(\{n\}-\{m\}\big)d_{n+m}.$$
    This bracket satisfies skew-symmetry
    $$\langle d_n,d_m\rangle_\s=-\langle d_m,d_n\rangle_\s,$$
and a $\s$-deformed Jacobi-identity
    \begin{align*}
        (q^n+1)\big \langle d_n, \langle d_l, d_m\rangle_\s\big
        \rangle_\s&+ (q^l+1)\big \langle d_l,  \langle d_m, d_n
\rangle_\s\big\rangle_\s+\notag\\
    &+(q^m+1)\big \langle d_m, \langle d_n,d_l\rangle_\s\big \rangle_\s\quad=0.
    \end{align*}
\end{thm}
\begin{remark}
    The linear space $\mfD_\sigma(\A)$ is a hom-Lie algebra with
    bilinear bracket defined in Theorem 3 and
    $\alpha:\mfD_\sigma(\A)\to\mfD_\sigma(A)$ given by
    $$\alpha(d_n)=\alpha(-t^nD)=\sigma(-t^n)D=-q^nt^nD=q^nd_n.$$
\end{remark}
\subsubsection{Non-linearly deformed Witt algebras}\label{sec:nonlinC}
The most general non-zero endomorphism $\s$ on $\k[t,t^{-1}]$ is
one on the form \mbox{$\sigma(t)=qt^s$} for $s\in\Z$ and
$q\in\k\setminus\{0\}$. With this $\sigma$, the left $\A$-module
$\mfD_\sigma(\A)$ is generated by a single element
$$D=\eta t^{-k+1}\frac{\id-\sigma}{t-qt^s}, \qquad \eta\in\k.$$ The element $\delta$
for this $D$ such that (\ref{eq:GenWittCondC2}) holds is
$$\delta=q^kt^{k(s-1)}\sum_{r=0}^{s-1}(qt^{s-1})^r.$$Equation
(\ref{eq:GenWittCondC1}) is clearly still valid. We then get,
using Theorem \ref{thm:GenWittC}, the following qhl-algebra
structure on $\mfD_\sigma(\A)$.
\begin{thm}\emph{(cf. \cite{HartLarsSilv1C})} \label{th:qWittnonlinearC} Let $\A=\k[t,t^{-1}]$. Then the $\k$-linear space
    $$\mfD_\s(\A) = \bigoplus_{n\in\Z} \k\cdot d_n,$$
    with
     $d_n=-t^nD$, can be
equipped
    with the bracket product
    $$\langle\,\cdot,\,\cdot\rangle_\s\,:\, \mfD_\s(\A)\times \mfD_\s(\A)\longrightarrow
    \mfD_\s(\A)$$ defined on generators by {\rm (\ref{eq:GenWittProdDefC})} as
    $$\langle d_n,d_m\rangle_\s=q^nd_{ns}d_m-q^md_{ms}d_n$$ and satisfying
    defining commutation relations
        \begin{align*}
            &\langle d_n,d_m\rangle _\sigma=\eta\sign(n-m)\sum_{l=\min(n,m)}^{\max(n,m)-1}q^{n+m-1-l}
            d_{s(n+m-1)-(k-1)-l(s-1)}\\
            &\qquad\text{for } n,m\geq 0;\\
            &\langle d_n,d_m\rangle_\sigma=\eta\Big
            (\sum_{l=0}^{-m-1}q^{n+m+l}d_{(m+l)(s-1)+ns+m-k}+\sum_{l=0}^{n-1}q^{m+l}d_{(s-1)l+n+ms-k}\Big )\\
            &\qquad\text{for } n\geq 0, m<0;\\
            &\langle d_n,d_m\rangle_\sigma=\eta\Big (\sum_{l_1=0}^{m-1}q^{n+l_1}
        d_{(s-1)l_1+m+ns-k}+
         \sum_{l_2=0}^{-n-1}q^{m+n+l_2}d_{(n+l_2)(s-1)+n+ms-k}\Big )\\
            &\qquad\text{for } m\geq 0, n<0;\\
            &\langle d_n,d_m\rangle _\sigma=\eta\sign(n-m)\sum_{l=\min(-n,-m)}^{\max(-n,-m)-1}q^{n+m+l}d_{(m+n)s+(s-1)l-k}\\
            &\qquad\text{for } n,m<0.
        \end{align*}
    Furthermore, this bracket satisfies skew-symmetry
    $$\langle d_n,d_m\rangle _\s=-\langle d_m,d_n\rangle_\s,$$ and a $\s$-deformed Jacobi-identity
    \begin{align*}
        \cs_{n,m,l}\,\Big (q^n\big \langle d_{ns},\langle d_m,d_l\rangle_\s\big
\rangle_\s+q^kt^{k(s-1)}\sum_{r=0}^{s-1}(qt^{s-1})^r\big \langle
        d_n,\langle d_m,d_l\rangle_\s\big \rangle_\s\Big )=0.
    \end{align*}
\end{thm}

\subsubsection{A generalization to several variables}

As shown in \cite{HartLarsSilv1C}, by not restricting the
denominator in $D\in\mfD_\sigma(\A)$ to \emph{greatest} common
divisors (on a certain subset of $\A$), but instead considering
only common divisors on this subset, we get, not the whole
$\A$-module $\mfD_\sigma(\A)$, but a proper submodule $\M$.
Theorem \ref{thm:GenWittC} can now be used to endow this submodule
with the structure of a qhl-algebra. For this, we need not
restrict our attention to one-variable Laurent polynomials.
\begin{thm}\emph{(cf. \cite{HartLarsSilv1C})} \label{th:qWittmultdimmonomC}
    Let $\A=\k[z_1^ {\pm 1},\dots z_n^{\pm 1}]$. The $\k$-linear space $$\M=
    \bigoplus_{\lbf\in\Z^n} \k\cdot d_{\lbf}$$
    spanned by $d_{\lbf}=-z_1^{l_1}\cdots z_n^{l_n}D$, where $D$ is given by
    $$D=Q\frac{\id-\s}{z_1^{G_1}\cdots z_n^{G_n}},$$ with $Q\in\k\setminus\{0\}$ and
    $(G_1,\dots, G_n)\in\Z^n$,
    can be endowed with a bracket
    defined on generators (by {\rm (\ref{eq:GenWittProdDefC})}) as
    \begin{align*}
        \langle d_{\kbf},d_{\lbf}\rangle_\s=q_{z_1}^{k_1}\cdots q_{z_n}^{k_n}d_{\a_1(\kbf),\dots,
\a_n(\kbf)}d_{\lbf}-
        q_{z_1}^{l_1}\cdots q_{z_n}^{l_n}d_{\a_1(\lbf),\dots, \a_n(\lbf)}d_{\kbf}
    \end{align*}
    and satisfying relations
    \begin{align*}
    \langle d_{\kbf},d_{\lbf}\rangle_\s&=
        Qq_{z_1}^{l_1}\cdots q_{z_n}^{l_n}d_{\alpha_1(\lbf)+k_1-G_1,\dots,
\alpha_n(\lbf)+k_n-G_n}-\\
        &\quad-Qq_{z_1}^{k_1}\cdots
        q_{z_n}^{k_n}d_{\alpha_1(\kbf)+l_1-G_1,\dots,\alpha_n(\kbf)+l_n-G_n}.
    \end{align*} The bracket satisfies the $\s$-deformed Jacobi-identity
    \begin{align*}
        \cs_{\kbf,\lbf,\hbf}\Big (
        &q_{z_1}^{k_1}\cdots q_{z_n}^{k_n}\big
\langle d_{\alpha_1(\kbf),\dots,\alpha_n(\kbf)},\langle
d_{\lbf},d_{\hbf}\rangle_\s \big
        \rangle_\s+\\
        &+q_{z_1}^{G_1}\cdots q_{z_n}^{G_n}z_1^{\d_1}\cdots z_n^{\d_n} \big
\langle d_{\kbf},\langle d_{\lbf},d_{\hbf}\rangle_\s \big
\rangle_\s
        \Big )=0.
    \end{align*}
\end{thm}
\subsection{A $(\alpha,\beta,\o)$-deformed loop algebra}
Loop algebras are special cases of Lie algebras known as current
algebras. Given a topological space $X$ and a Lie algebra $\g$ the
  \emph{current algebra with respect to $X$ and $\g$} consists of maps
  $X\to\g$ with (Lie) bracket multiplication
  \begin{align}\label{eq:loop_brack_class}
  [f,g](x)=[f(x),g(x)]_\g, \qquad \text{for } x\in X.
  \end{align}
  In the classical definition of current algebras it is usually
  required that $X$ is a smooth manifold, $\g$ is a
  finite-dimensional Lie algebra, and all maps belonging to the
  current algebra are smooth as well \cite{FuksC}. These smoothness assumptions
  and even the assumption that $X$ is a topological space and not
  just a set, are not important for the algebraic constructions in
  this article, since we do not discuss the topologies on the current
  algebra. It is enough for us that the set of maps comprising the
  algebra is a linear space closed under the commutator bracket
  multiplication (\ref{eq:loop_brack_class}).
   Taking $X$ to be the unit circle $S^1$, and restricting
to (trigonometric) polynomials we get the \emph{loop algebra of
$\g$}, oftentimes denoted by $\gh$. It is not difficult to see
that $\gh=\g\otimes\C[t,t^{-1}]$ with bilinear multiplication
determined by
$$[f\otimes t^n, g\otimes t^m]_{\gh}=[f,g]_\g\otimes t^{n+m}.$$
Loop algebras are important in physics, especially in conformal
field theories and superstring theory
\cite{QuantFieldStrings,DiFranMathiSenC,Fuchs1C,Fuchs2C}.

Let $\g$ be a qhl-algebra $(\g,\langle\cdot,\cdot \rangle_{\g},
\alpha_\g,\b_\g,\o_\g)$. Assume $\o_\g$ is linear over the tensor
product $\otimes:=\otimes_\k$ and form the vector space
$$\gh\overset{\text{def}}{=}\g\otimes \k[t,t^{-1}].$$ This means that
$\gh$ can be considered as the algebra of Laurent polynomials with
coefficients in the qhl-algebra $\g$. Define further
\begin{align*}
    \a_{\gh} &\overset{\text{def}}{=} \a_\g\otimes\id\\
    \b_{\gh} &\overset{\text{def}}{=} \b_\g\otimes\id\\
    \o_{\gh} &\overset{\text{def}}{=} \o_\g\otimes\id.
\end{align*}Furthermore, define a product on $\gh$ as
$$\langle x\otimes t^n, y\otimes t^m\rangle_{\gh}=\langle x,y\rangle_\g\otimes t^{n+m}.$$
\begin{thm}The vector space $\gh$ admits the structure of a qhl-algebra \newline\mbox{$(\gh,\langle\cdot,\cdot\rangle_{\gh},
\a_{\gh},\b_{\gh},\o_{\gh})$}.
\end{thm}
\begin{proof}The proof consists of a checking the axioms from Definition 1 of
qhl-algebras.
\begin{enumerate}
    \item The $\o_{\gh}$-skew symmetry is checked as follows:
    \begin{align*}
        \langle x\otimes t^n,y\otimes t^m\rangle_{\gh}&=\langle x,y\rangle_\g\otimes t^{n+m}=\big
        (\o_\g(x,y)\langle y,x\rangle_\g\big )\otimes t^{n+m}=\\
        &=\big (\o_\g(x,y)\otimes\id\big )\big (\langle
        y,x\rangle_\g\otimes t^{n+m}\big)=\\
&=\o_{\gh}(x,y)\langle y\otimes
        t^m,x\otimes t^n\rangle_{\gh}.
    \end{align*}
    \item Next we prove the $\b_{\gh}$-twisting of $\a_{\gh}$. First
          $\a_{\gh}(x\otimes t^n)=\a_\g(x)\otimes t^n$, and so
    \begin{align*}
        \langle \a_{\gh}(x\otimes t^n),\a_{\gh}(y\otimes t^m)\rangle_{\gh}&=\langle \a_\g(x)\otimes t^n,
        \a_\g(y)\otimes t^m\rangle_{\gh}=\\
        &=\langle \a_\g(x),\a_\g(y)\rangle_\g\otimes t^{n+m}=\\
        &=\big (\b_\g\circ\a_\g\langle x,y\rangle_\g\big
        )\otimes t^{n+m}=\\
        &=\b_{\gh}\circ\big (\a_\g\langle x,y\rangle_\g\otimes t^{n+m}\big
        )=\\
        &=\b_{\gh}\circ\a_{\gh}\big (\langle x,y\rangle_\g\otimes t^{n+m}\big ).
    \end{align*}
    \item The qhl-Jacobi identity lastly, is as follows. The
    left hand side is
    \begin{align*}
        \cs\, \o_{\gh}(z\otimes t^l,x\otimes t^n) &\Big (\big\langle \a_{\gh}(x\otimes t^n),\langle y\otimes
        t^m,z\otimes t^l\rangle_{\gh}\big \rangle_{\gh}+\\
        &+\b_{\gh}\big\langle x\otimes t^n,\langle y\otimes
        t^m,z\otimes t^l\rangle_{\gh}\big \rangle_{\gh}\Big ),
    \end{align*} where the
      notation $\cs$ here is used for cyclic summation
    with respect to $x\otimes t^n,y\otimes t^m,z\otimes t^l$. The first term in the parantheses is
    \begin{align*}
        \big\langle \a_{\gh}(x\otimes t^n),\langle y\otimes
        t^m,z\otimes t^l\rangle_{\gh}\big \rangle_{\gh}&=\big\langle \a_\g(x)\otimes t^n,\langle y\otimes
        t^m,z\otimes t^l\rangle_{\gh}\big \rangle_{\gh}=\\
        &=\big \langle \a_\g(x),\langle y,z\rangle_\g\big \rangle_\g\otimes t^{n+m+l}
    \end{align*} and the second
    \begin{align*}
        \b_{\gh}\big\langle x\otimes t^n,\langle y\otimes
        t^m,z\otimes t^l\rangle_{\gh}\big \rangle_{\gh}&=\b_{\gh} \big (\big \langle x,\langle
        y,z\rangle_\g\big \rangle_\g\otimes t^{n+m+l}\big )=\\
        &=\big (\b_\g\big \langle x,\langle y,z\rangle_\g\big
        \rangle_\g\big)\otimes t^{n+m+l}.
    \end{align*}Adding up these and using that $\o_{\gh}(z\otimes t^l,x\otimes t^n)=\o_\g(z,x)\otimes \id$,
    we get
    \begin{align*}
        &\cs\big(\o_\g(z,x)\otimes \id\big)\Big (\big \langle \a_\g(x),\langle y,z\rangle_\g\big \rangle_\g
        +\b_\g\big \langle x,\langle y,z\rangle_\g\big \rangle_\g\Big )\otimes t^{n+m+l}=\\
        &=\Big (\cs_{x,y,z}\,\o_\g(z,x)\big (\big \langle \a_\g(x),\langle y,z\rangle_\g\big \rangle_\g+
        \b_\g\big \langle x,\langle y,z\rangle_\g\big \rangle_\g\Big )\otimes t^{n+m+l}
    \end{align*} and the parantheses is zero since $\g$ is a qhl-algebra.
\end{enumerate}
\end{proof}

\section{Extensions}
Throughout this section we use that exact sequences of linear
spaces
\begin{align}\label{diag:vecExt1}
\begin{CD}
0 @>>> \afrak @>\iota>> E @>\pr>> L @>>> 0\\
\end{CD}
\end{align}\emph{split} in the sense that there is
a $\k$-linear map $s: L\to E$ called a \emph{section} such that
$\pr\circ s=\id_L$.
 Note that the condition $\pr\circ s=\id_L$ means that the $\k$-linear section $s$ is
one-to-one and so $L\cong s(L)$ (as linear spaces). This, together
with the exactness, lets us deduce that $E\cong
s(L)\oplus\i(\afrak)$ as linear spaces. Hence a basis of $E$ can
be chosen such that any $e\in E$ can be decomposed as
$e=s(l)+\i(a)$ for $a\in \afrak$ and $l\in L$, that is, we
consider $\i(\afrak)$ and $s(L)$ as subspaces of $E$.

To return to extensions, let us from now on assume that $L,\afrak$
and $E$ from (\ref{diag:vecExt1}) are qhl-algebras, and that we
have a section $s:L\to E$ such that $\o_E$ and $\o_L$ are
intertwined with $s$, that is
\begin{align}
    \label{eq:extcond2C}&\o_E\big(s(x)+\iota(a),s(y)+\iota(b)\big)\circ s=s\circ\o_L(x,y),
\end{align} if $(x,y)\in D_{\o_L}$ and $\big(s(x)+\iota(a),s(y)+\iota(b)\big)\in
D_{\o_E}$. In particular, $$\o_E(s(x),s(y))\circ
s=s\circ\o_L(x,y)$$ if $a$ and $b$ are taken to be zero. By
definition of a section $\pr\circ s=\id_L$, and since $\pr$ is a
homomorphism of algebras we have
$$0=\pr\circ (\langle s(x),s(y)\rangle_E-s\langle x,y\rangle_L)$$ which gives that
\begin{align}\label{eq:g_s}
    \langle s(x),s(y)\rangle_E=s\langle x,y\rangle_L+\iota\circ g(x,y),
\end{align} where $g:L\times L\to \afrak$ is a 2-cocycle-like $\k$-bilinear map which depends in a crucial way
on the section $s$. Thus $g$ is a measure of the deviation of $s$
from being a weak qhl-algebra morphism. Furthermore, $g$ has to
satisfy a generalized skew-symmetry condition on $D_{\o_L}$
\begin{align}\label{eq:g_e}
    \i\circ g(x,y)&=\langle s(x),s(y)\rangle_E-s\langle x,y\rangle_L=\notag\\
    &=\o_E(s(x),s(y))\circ \langle
    s(y),s(x)\rangle_E-s\circ\o_L(x,y)\circ \langle y,x\rangle_L=\notag\\
    &=\o_E(s(x),s(y))\circ\big (\langle s(y),s(x)\rangle_E-s\langle y,x\rangle_L
    \big )=\notag\\
    &=\o_E(s(x),s(y))\circ\i\circ g(y,x),
\end{align}for $(x,y)\in D_{\o_L}$ such that $(s(x),s(y))\in D_{\o_E}$.
Here we have used that $s$ intertwines $\o_E$ and $\o_L$.
\begin{dfn}
Denote the set of all maps $L\times L\to \afrak$ satisfying
{\rm(\ref{eq:g_e})} by $\Alt^2_\o(L,\afrak;\U)$, the
\emph{$\o$-alternating $\k$-bilinear
  mappings} associated with the extension {\rm(\ref{diag:vecExt1})}, which we denote by $\U$ to keep
  notation short.
\end{dfn}
\begin{remark}
    For Lie algebras $\o_E(s(x),s(y))$ is just multiplication by
    $-1$ and thus by linearity and injectivity of the map $\iota$,
    the condition {\rm(\ref{eq:g_e})} reduces to $$g(x,y)=-g(y,x),$$
    which is the classical skew-symmetry, independent on the
    extension.
\end{remark}

 By the commutativity of the boxes in (\ref{diag:ses_ab1}) we
have
$$\a_L\circ\pr =\pr\circ\a_E$$ which means that
$$\pr\circ (\a_E-s\circ\a_L\circ\pr)=0$$ and so
\begin{align}\label{eq:a_E}
    \a_E=s\circ\a_L\circ\pr+\i\circ f,
\end{align}where $f:E\to\afrak$ is a $\k$-linear map. By a similar argument we get
\begin{align}\label{eq:b_E}
    \b_E=s\circ\b_L\circ\pr+\i\circ h,
\end{align}for a $\k$-linear $h:E\to\afrak$. Obviously, both $f$ and $h$ depends on the section
chosen. To simplify notation we do not indicate explicitly this
dependence in what follows.

Since any $e\in E$ and $e'\in E$ can be decomposed as
$e=s(x)+\i(a)$ and $e'=s(y)+\iota(b)$ with $x,y\in L$ and
$a,b\in\afrak$, we have
\begin{align*}
    \langle e,e'\rangle_E&=\langle s(x)+\i(a),s(y)+\i(b)\rangle_E=\\
            &=\langle \i(a),\i(b)\rangle_E+\langle
            s(x),\i(b)\rangle_E+\langle\i(a),s(y)\rangle_E+\langle s(x),s(y)\rangle_E.
\end{align*}
With $\langle s(x),s(y)\rangle_E=s\langle x,y\rangle_L+\iota\circ
g(x,y)$ we can re-write this noting that by definition $\i$ is a
morphism of algebras:
\begin{align*}
        \langle e,e'\rangle_E&=\i\big( \langle
        a,b\rangle_\afrak\big )+\langle s(x),\i(b)\rangle_E+\langle
        \i(a),s(y)\rangle_E+s\langle x,y\rangle_L+\i\circ g(x,y)=\\
        &=s\langle x,y\rangle_L+\Big (\i\big ( \langle a,b\rangle
        _\afrak\big )+\langle s(x),\i(b)\rangle_E+\langle \i(a),s(y)\rangle_E+\i\circ g(x,y)\Big )
\end{align*} where the expression in parentheses is in $\i(\afrak)$ since $\i(\afrak)$ is an ideal in $E$
by the exactness. The extension is called \emph{inessential} if
$g\equiv 0$, which is equivalent to viewing $L$ as a subalgebra of
$E$.

We consider only \emph{central} extensions, i.e. extensions
satisfying
$$\iota(\afrak)\subseteq Z(E)=\{e\in E\, |\,
\langle e,E\rangle_E=0\,\},$$ where $\afrak$ is \emph{abelian},
that is $\langle\afrak,\afrak\rangle_\afrak=0.$ This means in
particular that
\begin{align*}
        \langle e,e'\rangle_E&=s\langle x,y\rangle_L+\Big
        (\underbrace{\i\big( \langle
          a,b\rangle_\afrak\big)}_{=0}+\underbrace{\langle
          s(x),\i(b)\rangle_E}_{=0}+
        \underbrace{\langle\i(a),s(y)\rangle_E}_{=0}+\i\circ g(x,y)\Big )=\\
        &=s\langle x,y\rangle_L+\i\circ g(x,y).
\end{align*}
 To study the necessary conditions on $g$, when having a Jacobi-like
 identity on $E$,
 under the centrality assumption, we calculate
\begin{multline*}
    \o_E\big(s(z)+\iota(c),s(x)+\iota(a)\big)\Big (\big\langle \a_E(s(x)),\langle
s(y),s(z)\rangle_E\big \rangle_E+\\
+\b_E \big \langle s(x),\langle s(y),s(z)\rangle_E\big
\rangle_E\Big )
\end{multline*} and sum up this cyclically.
The reason for not considering products of elements on the form
$s(x)+\i(a)$ in the brackets is that, since the extension is
assumed to be central, the terms involving elements from
$\iota(\afrak)$ within the brackets vanish automatically. So, we
have
\begin{multline*}
    \big\langle \a_E(s(x)),\langle s(y),s(z)\rangle_E\big
    \rangle_E=\big\langle \a_E(s(x)),s\langle y,z\rangle_L+\i\circ
    g(y,z)
\big \rangle_E=\\
    =\big\langle s(\a_L(x))+\i\circ f(s(x)),s\langle y,z\rangle_L+\i\circ g(y,z)\big\rangle_E=\\
    =\big\langle s(\a_L(x)),s\langle y,z\rangle_L\big\rangle_E=\\
    =s\big \langle\a_L(x),\langle y,z\rangle_L\big \rangle_L+
\i\circ g(\a_L(x),\langle y,z\rangle_L).
\end{multline*}In a similar fashion we see that
\begin{align*}
    \b_E\big \langle s(x),\langle
    s(y),s(z)\rangle_E\big\rangle_E&=
\b_E\big \langle s(x),s\langle y,z\rangle_L+\i\circ g(y,z)\big
\rangle_E
=\\
    &=\b_E\circ s\big \langle x,\langle y,z\rangle_L\big
    \rangle_L+\b_E\circ\i\circ g(x,\langle y,z\rangle_L).
\end{align*}Noting that by exactness
$$\b_E\circ\i\circ g=\i\circ h\circ\i\circ g,$$
we can re-write the above as
\begin{align*}
    &\b_E\circ s\big \langle x,\langle y,z\rangle_L\big
    \rangle_L+\b_E\circ\i\circ g(x,\langle y,z\rangle_L)=\\
    &=s\circ\b_L\big \langle x,\langle y,z\rangle_L\big
    \rangle_L+\i\circ h\circ s\big \langle x,\langle y,z\rangle_L\big
    \rangle_L+\i\circ h\circ\i\circ g(x,\langle y,z\rangle_L)=\\
    &=s\circ\b_L\big \langle x,\langle y,z\rangle_L\big
    \rangle_L+\i\circ h\circ\big (s\big \langle x,\langle
    y,z\rangle_L\big \rangle_L+\i\circ g(x,\langle y,z\rangle_L)\big
    )
\end{align*}and so the qhl-Jacobi identity becomes
\begin{multline*}
    0=\cs_{(x,a),(y,a),(z,c)}\,\o_E\big(s(z)+\iota(c),s(x)+\iota(a)\big)\circ\\
    \circ\Big (s\big \langle\a_L(x),\langle
    y,z\rangle_L\big \rangle_L+\i\circ g(\a_L(x),\langle y,z\rangle_L)+\\
    +s\circ\b_L\big \langle x,\langle
    y,z\rangle_L\big\rangle_L+\i\circ h\circ\big (s\big \langle
    x,\langle y,z\rangle_L\big \rangle_L+\i\circ g(x,\langle y,z\rangle_L)\big
    )\Big )=\\
    =\cs_{(x,a),(y,a),(z,c)}\,\o_E\big(s(z)+\iota(c),s(x)+\iota(a)\big)\circ\Big (\i\circ g(\a_L(x),\langle y,z\rangle_L)+\\
    +\i\circ h\circ\big (s\big \langle x,\langle y,z\rangle_L\big
    \rangle_L+\i\circ g(x,\langle y,z\rangle_L)\big
    )\Big )
\end{multline*} since $L$ is a qhl-algebra. Note that $s\big \langle
x,\langle y,z\rangle_L\big \rangle_L+\i\circ g(x,\langle
y,z\rangle_L)$ can be re-written as $\big \langle s(x),s\langle
y,z\rangle_E\big \rangle_E,$ which will be useful later. So, a
necessary condition for $E$ to be a qhl-algebra is that
\begin{align}\label{eq:jacobi_gC}
    \cs_{(x,a),(y,a),(z,c)}\,&\o_E\big(s(z)+\iota(c),s(x)+\iota(a)\big)\circ\Big (\i\circ g(\a_L(x),\langle
    y,z\rangle_L)+\notag\\
    &+\i\circ h\circ\big (s\big \langle x,\langle y,z\rangle_L\big
    \rangle_L
+\i\circ g(x,\langle y,z\rangle_L)\big
    )\Big )=0.
\end{align}Or, after re-writing the second term,
\begin{multline}\label{eq:rewrite_gC}
\cs_{(x,a),(y,a),(z,c)}\,\o_E\big(s(z)+\iota(c),s(x)+\iota(a)\big)\circ\\
\circ\Big (\i\circ g(\a_L(x),\langle y,z\rangle_L)+\i\circ
h\circ\big \langle s(x),s\langle y,z\rangle_E\big \rangle_E
    \Big )=0.
\end{multline} We denote the subset of all $\o$-alternating $g\in\Alt^2_\o(L,\afrak;\U)$,
$g:L\times L\to \afrak$ satisfying the condition
    (\ref{eq:jacobi_gC}) by $Z^2_{\o,h}(L,\afrak;\U)$, the set of
    \emph{2-cocycle-like maps}, by analogy with the Lie algebra case.
    The subscript $h$ is to remind us that the 2-cocycle-like maps may depend on $h$.
    We will now show that (\ref{eq:jacobi_gC}), or equivalently
(\ref{eq:rewrite_gC}), is independent of the choice of section $s$
and the $h$.

Taking another section $\tilde s$ with $\pr\circ \tilde s =\id_L$
satisfying the intertwining condition $\o_E\big (\tilde
s(x)+\iota(a),\tilde s(y)+\iota(b)\big)\circ\tilde s=\tilde
s\circ\o_L(x,y)$, we see that
$$(\tilde s-s)(x)=\i\circ k(x),$$ for some linear $k:L\to \afrak$, and so
\begin{align}\label{eq:skC}
    \tilde s=s+\i\circ k.
\end{align}Hence
\begin{align*}
    \iota\circ\tilde g(x,y)&=\langle \tilde s(x),\tilde s(y)\rangle_E-\tilde
    s\langle x,y\rangle_L =\langle s(x)+\i\circ k(x),s(y)+\i\circ
    k(y)\rangle_E -\\
    &-s\langle x,y\rangle_L -\i\circ k(\langle x,y\rangle
    )=\iota\circ g(x,y)-\i\circ k(\langle x,y\rangle_L )
\end{align*} since the extension is central. By the injectivity of
$\iota$ we get $$\tilde g(x,y)=g(x,y)-k(\langle x,y\rangle_L).$$
Further
\begin{align*}
    \iota\circ(\tilde h-h)(x)&=\big (\b_E-\tilde s\circ \b_L\circ\pr-\b_E+s\circ\b_L\circ\pr\big )(x)=\\
    &=(s-\tilde s)\circ \b_L\circ\pr(x)=-\i\circ k\circ\b_L\circ\pr(x)
\end{align*} giving since $\iota$ is an injection
\begin{align}\label{eq:tildew_wC}
    \tilde h =h-k\circ\b_L\circ\pr.
\end{align}Therefore
\begin{multline*}
    \cs_{(x,a),(y,a),(z,c)}\Big (\tilde g(\a_L(x),\langle y,z\rangle_L)+\tilde
    h\big \langle \tilde s(x),\langle \tilde s(y),\tilde s(z)\rangle_E\big \rangle_E\Big )=\\
    =\cs_{(x,a),(y,a),(z,c)}\Big (g(\a_L(x),\langle y,z\rangle_L)-
    k\big \langle \a_L(x),\langle y,z\rangle_L\big \rangle_L+\\
     +h\big \langle s(x),\langle s(y),s(z)\rangle_E\big \rangle_E -k\circ\b_L\big
    \langle x,\langle y,z\rangle_L\big \rangle_L\Big )=\\
    =\cs_{(x,a),(y,a),(z,c)}\Big (g(\a_L(x),\langle y,z\rangle_L)+h\big
    \langle s(x),\langle s(y),s(z)\rangle_E\big \rangle_E\Big )
\end{multline*}since $L$ is a qhl-algebra.
Knowing that
    $$\langle \a_E(x),\a_E(y)\rangle_E=\b_E\circ\a_E\langle x,y\rangle_E$$ by the definition of a qhl-algebra, we can deduce another
relation that necessarily must hold for $E$ to be a qhl-algebra.
The left hand side of the above equality for $x=s(x)$ and $y=s(y)$
can be written as follows:
\begin{multline*}
    \langle\a_E(s(x)),\a_E(s(y))\rangle_E= \langle
    (s\circ\a_L\circ\pr+\i\circ f)\circ
    s(x),(s\circ\a_L\circ\pr+\i\circ f)\circ s(y)\rangle_E=\\
        =\langle s\circ\a_L(x),s\circ\a_L(y)\rangle_E=s\langle \a_L(x),\a_L(y)\rangle_L+\\
        +\i\circ g(\a_L(x),\a_L(y)),
\end{multline*} where we have used the centrality in the second equality. The right hand side is
\begin{align*}
    \b_E\circ\a_E\langle s(x),s(y)\rangle_E&=\b_E\circ
    s\circ\a_L\langle x,y\rangle_L+\b_E\circ\i\circ
    f\langle s(x),s(y)\rangle_E=\\
        &=(s\circ\b_L\circ\pr+\i\circ h)\circ s\circ\a_L\langle
        x,y\rangle_L+\\
        &\quad+(s\circ\b_L\circ\pr+\i\circ h)\circ\i
        \circ f\langle s(x),s(y)\rangle_E=\\
        &=s\circ\b_L\circ\a_L\langle x,y\rangle_L+\i\circ h\circ
        s\circ\a_L\langle x,y\rangle_L+\\
        &\quad +\i\circ h\circ\i\circ f\langle s(x),s(y)\rangle_E
\end{align*}and after comparing and using injectivity of $\iota$, we get
that
\begin{align*}
    g(\a_L(x),\a_L(y))=h\circ \big (s\circ\a_L\langle
    x,y\rangle_L+\i\circ f\langle s(x),s(y)\rangle_E\big ).
\end{align*}Hence, we get the following result.
\begin{thm}\label{thm:qhlext1C}
    Suppose $(L,\a_L,\b_L,\o_L)$ and $(\afrak,\a_\afrak,\b_\afrak,\o_\afrak)$ are qhl-algebras with $\afrak$ abelian. If there exists a central extension
    $(E,\a_E,\b_E,\o_E)$ of $(L,\a_L,\b_L,\o_L)$ by
    $(\afrak,\a_\afrak,\b_\afrak,\o_\afrak)$ then for a section $s:L\to E$, satisfying
    {\rm (\ref{eq:extcond2C})}
    there is an \mbox{$\o$-alternating} bilinear $g:L\times L\to\afrak$ and
    linear maps $f:E\to\afrak$ and $h:E\to\afrak$ such that
    \begin{align}
        &f\circ\iota=\a_\afrak \label{eq:fanda_aC},\\
        &h\circ\iota=\b_\afrak \label{eq:handb_aC}
    \end{align}and also
    \begin{align}\label{eq:both_g_and_f_and_hC}
        g(\a_L(x),\a_L(y))=h\circ \big (s\circ\a_L\langle
        x,y\rangle_L+\i\circ f\langle s(x),s(y)\rangle_E\big )
    \end{align}
and
\begin{align}\label{eq:twococycle2C}
    \cs_{(x,a),(y,a),(z,c)}\,\o_E&\big(s(z)+\iota(c),s(x)+\iota(a)\big)\circ \Big (\i\circ g(\a_L(x),\langle
    y,z\rangle_L)+\notag\\
    &+\i\circ h\circ\big (s\big \langle x,\langle y,z\rangle_L\big
    \rangle_L+\i\circ g(x,\langle y,z\rangle_L)\big
    )\Big )=0
\end{align} for all pairs ${(x,a),(y,b),(z,c)}\in L\times\afrak$ such
that $\big(s(z)+\iota(c),s(x)+\iota(a)\big)$,
$\big(s(x)+\iota(a),s(y)+\iota(b)\big),\big(s(y)+\iota(b),s(z)+\iota(c)\big)\in
D_{\o_E}$. Moreover, equation {\rm(\ref{eq:twococycle2C})} is
independent of the choice of section $s$ and function $h$, if only
sections $s$ satisfying {\rm (\ref{eq:extcond2C})} are considered.
\end{thm}
\begin{proof}
    It remains to prove (\ref{eq:fanda_aC}) and (\ref{eq:handb_aC}). However these follows directly from
    (\ref{eq:a_E}), (\ref{eq:b_E}), the commutativity of (\ref{diag:ses_ab1}) and injectivity
    of $\iota$.
\end{proof}
\begin{example}\label{ex:homLieExtC}
    By taking $\b_L=\id_L$, $\b_E=\id_E$, $\b_\afrak=\id_\afrak$ and $\o_L(x,y)v_L=-1\cdot v_L$ for all
    $x,y,v_L\in L$, $\o_E(e,e')v_E=-1\cdot v_E$ for all $e,e',v_E\in E$ (we have $D_{\o_L}=L\times L$ and $D_{\o_E}=E\times E$ here),
    that is if we consider only \mbox{hom-Lie algebras}, we recover the results from
    the previous paper {\rm\cite{HartLarsSilv1C}}. To see this consider first {\rm(\ref{eq:twococycle2C})}.
    The assumption that $\b_E=\id_E$ and $\b_L=\id_L$ implies that
    \begin{align}\label{eq:hid}
          \i\circ h=\id_E-s\circ\pr,
    \end{align} and hence by exactness
    \begin{align*}
        &\i\circ h\circ\big (s\big \langle x,\langle y,z\rangle_L\big
        \rangle_L+\i\circ g(x,\langle y,z\rangle_L)\big
         )=\\
        &=(\id_E-s\circ\pr)\circ\big (s\big \langle x,\langle
        y,z\rangle_L\big \rangle_L+\i\circ g(x,\langle y,z\rangle_L)\big
        )=\\
        &=s\big \langle x,\langle y,z\rangle_L\big\rangle_L-s\big
        \langle x,\langle y,z\rangle_L\big\rangle_L+\i\circ
        g(x,\langle y,z\rangle_L)=\\
        &=\iota\circ g(x,\langle y,z\rangle_L).
    \end{align*}This means that {\rm(\ref{eq:twococycle2C})} can be re-written
    using that $\iota$ is an injective qhl-algebra morphism as
    $$\cs_{x,y,z}\, g\big ((\id_L+\a_L)(x),\langle y,z\rangle_L\big
    )=0,$$ which is the condition obtained in the previous
    work {\rm\cite{HartLarsSilv1C}}. In the same manner one can see,
    using {\rm(\ref{eq:hid})} and injectivity of $\iota$ that {\rm(\ref{eq:both_g_and_f_and_hC})} reduces to
    $$g(\a_L(x),\a_L(y))= f\langle s(x),s(y)\rangle_E$$
which also is a formula from our previous article.

Note that {\rm(\ref{eq:hid})} can be written as
$h\circ\iota|_{\afrak}=\id_\afrak$ and $h\circ s|_{L}=0$. Indeed,
we can decompose any $e\in E$ as $e=s(x)+\iota(a)$ and so
\begin{align*}
     \iota\circ h(s(x)+\iota(a))&=(\id_E-s\circ\pr)(s(x)+\iota(a))=\\
          &=0+\iota(a)=\iota(a).
\end{align*} Since $\iota$ is an injection this gives
$h(s(x)+\iota(a))=\iota(a)$. Restricting even further to (color)
Lie algebras and thus having $\alpha_L=\alpha_E=\id$, we have $f$
satisfies a similar condition $f\circ\iota|_{\afrak}=\id_\afrak$
and $f\circ s|_{L}=0$. \qed
\end{example}
\begin{example}\label{ex:colorLieExt}
Consider the following short exact sequence of color Lie algebras
$L,E$ and $\afrak$ with the same $\Gamma$-grading and commutation
factor $\varepsilon$
$$\begin{CD}
0@>>>\afrak @>\iota>>E@>\pr >>L @>>> 0
\end{CD}$$with $\iota(\afrak)$ central in $E$. This setup is a special case
of the construction of Scheunert and Zhang
{\rm\cite{ScheunertZhangC}} and Scheunert {\rm\cite{Scheunert1C}},
special in the sense that we consider \emph{central} extensions
and not just \emph{abelian}. We shall show that our construction
encompasses the ones in {\rm\cite{Scheunert1C,ScheunertZhangC}}
for central extensions.

We first briefly recall Scheunert and Zhang's construction (as
given in {\rm\cite{Scheunert1C}}). Thus we consider our algebras
graded as
$$\afrak=\bigoplus_{\gamma\in\Gamma}\afrak_\gamma,\qquad L=\bigoplus_{\gamma\in\Gamma}L_\gamma,\qquad
    E=\bigoplus_{\gamma\in\Gamma}E_\gamma $$ and so the above sequence becomes
$$\begin{CD}
0@>>>\bigoplus_{\gamma\in\Gamma}\afrak_\gamma
@>\iota>>\bigoplus_{\gamma\in\Gamma}E_\gamma@>\pr
>>\bigoplus_{\gamma\in\Gamma}L_\gamma @>>> 0.
\end{CD}$$Note that this means that $\iota$ and $\pr$ are color Lie algebra homomorphisms and this in turn
implies that they are homogeneous of degree zero. Take a section
$s:L\to E$ which is \emph{homogeneous of degree zero}, that is, $s
(L_\gamma)\subseteq E_\gamma$ for all $\gamma\in\Gamma$. With this
data the Scheunert 2-cocycle condition can be expressed as
\begin{align}\label{eq:Sch-2-cocycle}
\cs_{x,y,z}\, \e(\gamma_z,\gamma_x)g(x,\langle y,z\rangle)=0,
\end{align} for homogenous elements $x,y,z$ and where $\gamma_x,\gamma_y,\gamma_z$ are the graded degrees of
$x,y,z$ respectively.

Putting the above in a qhl-algebras setting means letting $\o$
play the role of the commutation factor $\varepsilon$, where $\o$
is then defined on homogeneous elements,
\mbox{$D_\o=D_\varepsilon=\cup_{\gamma\in\Gamma}L_\gamma$,} and
dependent only on the graded degree of these elements.
 The set $\Alt^2_\varepsilon(L,\afrak;\U)$ includes all $g$ coming
 from the "defect"-relation $\iota\circ g(x,y)=\langle s(x),s(y)\rangle_E-s\langle
 x,y\rangle_L$, for $s$ a homogeneous section of degree zero. Hence all
 such $g$'s are also homogeneous of degree zero. Noting that
 $h\circ\iota|_{\afrak}=\id_\afrak$ and $h\circ s|_{L}=0$ from
    the Example \ref{ex:homLieExtC}, the relation
    {\rm(\ref{eq:twococycle2C})} now becomes,
    \begin{multline*}
       \cs_{(x,a),(y,a),(z,c)}\,\varepsilon_L\big(s(z)+\iota(c),s(x)+\iota(a)\big)\circ\\
       \circ\Big (\iota\circ g(\a_L(x),\langle
    y,z\rangle_L)+\iota\circ h\circ\big (s\big \langle x,\langle y,z\rangle_L\big
    \rangle_L+\i\circ g(x,\langle y,z\rangle_L)\big
    )\Big )=\\
       =2\cs_{(x,a),(y,a),(z,c)}\,\varepsilon\big(s(z)+\iota(c),s(x)+\iota(a)\big)\circ\iota\circ g(x,\langle
       y,z\rangle_L)=0
     \end{multline*} implying that
    $$\cs_{x,y,z}\,\varepsilon(z,x) g(x,\langle y,z\rangle_L)=0,$$ for homogeneous elements, which is
    Scheunert's 2-cocycle condition for central extensions {\rm\cite{Scheunert1C}}.
 \qed
\end{example}

\subsection{Equivalence between extensions}
Let $\varphi:E\to E'$ be a weak qhl-algebra morphism such that $E$
and $E'$ are of the same $\varphi$-related $\o$-type. We call two
extensions
\begin{align*}
\begin{CD}
    0 @>>> (\afrak,\Mor(\afrak)) @>\iota>> (E,\Mor(E)) @>\pr>> (L,\Mor(L)) @>>> 0
\end{CD}
\end{align*} and
\begin{align*}
\begin{CD}
    0 @>>> (\afrak,\Mor(\afrak)) @>\iota'>> (E',\Mor(E')) @>\pr'>> (L,\Mor(L)) @>>> 0
\end{CD}
\end{align*} \emph{weakly equivalent} or a \emph{weak equivalence} if the diagram
\begin{align}\label{diag:equiv}
\begin{CD}
    0 @>>> (\afrak,\Mor(\afrak)) @>\iota>> (E,\Mor(E)) @>\pr>> (L,\Mor(L)) @>>> 0\\
    @. @| @V\f VV  @| \\
    0 @>>> (\afrak,\Mor(\afrak)) @>\iota'>> (E',\Mor(E')) @>\pr'>> (L,\Mor(L)) @>>> 0\\
\end{CD}
\end{align} commutes. Similarly one defines \emph{strong equivalence} as a diagram with the map $E\to E'$
being a strong morphism. That $\f$ is automatically an isomorphism
of linear spaces follows from the $5$-lemma.
\begin{dfn}The set of weak
equivalence classes of extensions of $L$ by $\afrak$ is denoted by
$E^-(L,\afrak)$, and the set of strong equivalence classes by
$E^+(L,\afrak)$.
\end{dfn}
\begin{remark} In the case of Lie algebras, or generally, color Lie
  algebras, weak and strong extensions coincide since weak and strong
  morphisms do.
\end{remark}

We pick sections $s:L\to E$ and $s':L\to E'$ satisfying $\pr\circ
s=\id_L=\pr'\circ s'$ such that
$$\o_{E'}\big(s'(x)+\iota'(a),s'(y)+\iota'(b)\big)\circ
s'=s'\circ\o_L(x,y)$$ and
$$\o_E\big(s(x)+\iota(a),s(y)+\iota(b)\big)\circ
s=s\circ\o_L(x,y).$$ Then there is a $g'\in
\Alt_{\o_L}^2(L,\afrak;\U')$ associated with the extension $\U'$
of $L$ by $\afrak$ such that $\langle
s'(x),s'(y)\rangle_{E'}=s'\langle s,y\rangle_L+\iota'\circ
g'(x,y)$.

Given a map $\f: E\to E'$ such that the diagram (\ref{diag:equiv})
commutes means in particular that
$$\pr'\circ\f=\id_L\circ \pr$$ and so
$$\pr'\circ\f(s(x))=x,$$ which gives
$$0=\pr'\circ s'(x)-\pr'\circ\f(s(x))=\pr'\big (s'(x)-\f(s(x))\big ).$$ Hence
\begin{align}\label{eq:s's}
    s'(x)=\f\circ s(x)+\i'\circ\xi(x)
\end{align} for some $\k$-linear $\xi:L\to \afrak$. Taking $x,y\in L$ we have, using the centrality,
\begin{multline*}
    \i'\circ g'(x,y)=\langle s'(x),s'(y)\rangle_{E'}-s'\langle
    x,y\rangle_{E'}=\\
    =\langle\f(s(x))+\i'\circ\xi(x),\f(s(y))+\i'\circ\xi(y)\rangle_{E'}-
    \f( s\langle x,y\rangle_L)-\i'\circ\xi(\langle x,y\rangle_L)=\\
    =\f\circ\big (\langle s(x),s(y)\rangle_E-s\langle
    x,y\rangle_L\big )-\i'\circ\xi(\langle
    x,y\rangle_L)=\f\circ\i\circ g(x,y)-\i'\circ\xi(\langle x,y\rangle_L)
\end{multline*}and since $\f \circ \i=\i'$ by (\ref{diag:equiv}) we get
$$
   \i'\circ g'(x,y)=\i'\circ g(x,y)-\i'\circ\xi(\langle x,y\rangle_L)
$$ or
\begin{align}\label{eq:g'g}
      g'(x,y)=g(x,y)-\xi(\langle x,y\rangle_L),
\end{align}by the injectivity of $\i'$.

\begin{remark}
    Note that {\rm(\ref{eq:g'g})} implies
    \begin{multline*}
        (\o_{E'}(s'(y),s'(x))-\id_{E'})\circ\iota'\circ
        g'(x,y)=\\
        =(\o_E(\varphi\circ s(y),\varphi\circ s(x))-\id_{E'})\circ\iota'\circ
        g(x,y)-\iota'\circ\xi\circ(\o_L(y,x)-\id_L)(\langle x,y\rangle_L)
    \end{multline*} by exchanging $x$ and $y$ and subtracting
    equalities. For Lie algebras we get
    \begin{multline*}
        (-\id_{E'}-\id_{E'})\circ\iota'\circ
        g'(x,y)=(-\id_{E'}-\id_{E'})\circ\iota'\circ g(x,y)-\\
        -\iota'\circ\xi\circ(-\id_{L}-\id_{L})(\langle
        x,y\rangle_L).
    \end{multline*}This is obviously equivalent to
    $$\iota'\circ g'(x,y)=\iota'\circ g(x,y)-\iota'\circ\xi
    (\langle x,y\rangle_L)$$ or, by the injectivity of $\iota'$
    $$g'(x,y)=g(x,y)-\xi
    (\langle x,y\rangle_L)$$ which is nothing but the standard
    equivalence 2-cocycle-condition for the second cohomology group $H^2$ for Lie algebras.

    For color Lie algebra extensions we have the same commutation factor and $\Gamma$-grading
    for both $E,E'$ and $L,\afrak$. Hence
    \begin{align*}
        &\o_{E'}(s'(y),s'(x))v_{E'}=\varepsilon\big(\deg(s'(y)),\deg(s'(x))\big)v_{E'},\\
        &\o_{E'}(\varphi\circ s(y),\varphi\circ s(x))v_{E'}=
        \varepsilon\big(\deg(\varphi\circ s(y)),\deg(\varphi\circ s(x))\big)v_{E'},\\
        &\o_L(y,x)v_L=\varepsilon(\deg(y),\deg(x))v_L,
    \end{align*}on homogeneous elements, where $\deg(x)$ denotes the graded degree of $x$ and
    $v_L\in L, v_{E'}\in E'$. This means
    that
    \begin{multline*}
        \big
        (\varepsilon\big(\deg(s'(y)),\deg(s'(x))\big)-\id_{E'}\big
        )\circ \iota'\circ g'(x,y)=\\
        =\big (\varepsilon\big(\deg(\varphi\circ s(y)),\deg(\varphi\circ s(x))\big)
        -\id_{E'}\big )\circ\iota'\circ
        g(x,y)-\\
        -\iota'\circ\xi\big
        (\varepsilon(\deg(y),\deg(x))-\id_L\big )(\langle
        x,y\rangle_L)
    \end{multline*} for homogeneous elements $x,y$ and $z$.
    Since $s$ and $\varphi$ is of homogeneous degree zero this
    implies that
    \begin{align*}
    \varepsilon(x,y)&=\varepsilon\big(\deg(s'(y)),\deg(s'(x))\big)=
    \varepsilon\big(\deg(\varphi\circ s(y)),\deg(\varphi\circ
    s(x))\big)=\\
    &=\varepsilon(\deg(y),\deg(x))
    \end{align*} and so
    $$g'(x,y)=g(x,y)-\xi(\langle x,y\rangle_L)$$ for color Lie
    algebras also.
    \end{remark}
    We can view $\xi(\langle x,y\rangle)_L$ as a ''2-coboundary'' thus
motivating the following definition.
\begin{dfn}
The set of all 2-cocycle-like maps modulo  2-coboundary-like with
respect to a weak isomorphism is denoted by
$H^2_{\o,-}(L,\afrak;\U)$ and with respect to a strong isomorphism
$H^2_{\o,+}(L,\afrak;\U)$.
\end{dfn}

Now, given the relation $g'(x,y)=g(x,y)-\xi(\langle x,y\rangle_L)$
and two extensions $(E,\Mor(E))$ and $(E',\Mor(E'))$ of $L$ by
$\afrak$, can we construct a weak (strong) equivalence, that is, a
weak (strong) isomorphism making (\ref{diag:equiv}) commute? We
can view $E$ and $E'$ as
\begin{align*}
       &E= s(L)\oplus\i(\afrak)\qquad (\text{as linear spaces})\\
       &E'= s'(L)\oplus\i'(\afrak)\qquad (\text{as linear spaces})
\end{align*}
since the sequences on the form (\ref{diag:vecExt1}) is a split
sequence. This means that any element $e\in E$ can be decomposed
as
$$e=s(l)+\i(a)$$ for $a\in\afrak$ and $l\in L$. We define a map $\varphi:E\to E'$ by
\begin{align}\label{eq:f_def}
    \f\big (s(l)+\i(a)\big )\overset{\text{def}}{=}s'(l)+\i'(a-\xi(l)).
\end{align}We will show that this is a weak isomorphism of qhl-algebras. That it is surjective is clear. Suppose
that
    $$s'(l)+\i'(a-\xi(l))=s'(\tilde l)+\i'(\tilde a-\xi(\tilde l)).$$ This is equivalent to
    $$s'(l-\tilde l)+\i'(a-\tilde a+\xi(l-\tilde l))=0$$ and so injectivity follows from the injectivity of
    $\i'$ and $s'$.
To have a weak equivalence we must check
\begin{align}\label{eq:check1}
    \langle \f(x),\f(y)\rangle_{E'}=\f\langle x,y\rangle_E.
\end{align}
On the one hand,
\begin{multline*}
    \langle \f(s(l)+\i(a)),\f(s(\tilde l)+\i(\tilde
    a))\rangle_{E'}=\langle s'(l)+\i'(a-\xi(l)),s'(\tilde
    l)+\i'(\tilde a-\xi(\tilde l))\rangle_{E'}=\\
    =\langle s'(l),s'(\tilde l)\rangle_{E'}=s'\langle l,\tilde
    l\rangle_L+\i'\circ g'(l,\tilde l)=\\
    =s'\langle l,\tilde l\rangle_L+\i'\circ\big (g(l,\tilde
    l)-\xi\langle l,\tilde l\rangle_L\big ).
\end{multline*}On the other hand,
\begin{align*}
    \f\circ \langle s(l)+\i(a),s(\tilde l)+\i(\tilde
    a)\rangle_E&=\f\big (\langle s(l),s(\tilde l)\rangle_E\big )=
    \f\big (s\langle l,\tilde l\rangle_L +\i\circ g(l,\tilde l)\big )=\\
    &=s'\langle l,\tilde l\rangle_L+\i'\circ\big (g(l,\tilde
    l)-\xi\langle l,\tilde l\rangle_L\big )
\end{align*}from which we see that (\ref{eq:check1}) is proved. Hence,
\begin{thm} With definitions and notations as above, there is a one-to-one correspondence between elements of $E^-(L,\afrak)$ and elements of $H^2_{\o,-}(L,\afrak)$.
\end{thm}Rephrased, the theorem says that there is a one-to-one correspondence between weak equivalence classes of
central extensions of $L$ by an abelian $\afrak$ and bilinear maps
$g$ transforming according to (\ref{eq:g'g}) under weak
isomorphisms $E\to E'$.

If we are seeking strong equivalence we also have to condition and
check the intertwining conditions
\begin{description}
    \item[M2a.] $\f\circ\a_E=\a_{E'}\circ\f$, and
    \item[M3a.] $\f\circ\b_E=\b_{E'}\circ\f$.
\end{description} The check of M3a is completely analogous to M2a so we can restrict our attention to M2a.
To find conditions for M2a we first note that
$$f\circ\i=\a_{\afrak}=f'\circ\i'$$ which can be seen from the
commutativity of (\ref{diag:equiv}). Hence,
\begin{align}\label{eq:M2_1}
    \f\circ\a_E\big (s(l)+\i(a)\big )&=\f\big (s\circ\a_L(l)+\i\circ f(\i(a)+s(l))\big )=\notag\\
    &=s'\circ\a_L(l)+\i'\circ\big (f(\i(a)+s(l))-\xi(\a_L(l))\big )=\notag\\
    &=s'\circ\a_L(l)+\i'\circ\big (\a_{\afrak}(a)+f\circ s(l)-\xi(\a_L(l))\big ),
\end{align}and
\begin{align}\label{eq:M2_2}
    \a_{E'}\circ\f\big (s(l)+\i(a)\big )&=\a_{E'}\circ\big (s'(l)+\i'(a-\xi(l))\big )=\notag\\
    &=s'\circ\a_L(l)+\i'\circ f'\circ \big (s'(l)+\i'(a-\xi(l))\big )=\notag\\
    &=s'\circ\a_L(l)+\i'\circ\big (\a_{\afrak}(a)-\a_{\afrak}(\xi(l))+f'\circ s'(l)\big ).
\end{align} Combining
(\ref{eq:M2_1}) and (\ref{eq:M2_2}) shows that a sufficient
condition for M2a to hold is that
\begin{description}
    \item[S1$\a$.] $\a_{\afrak}\circ\xi=\xi\circ\a_{L}$ and
    \item[S2$\a$.] $f\circ s=f'\circ s'$
\end{description}or in commutative diagram form
$$\begin{CD}
L @>\xi>> \afrak \\
@V\a_{L} VV @VV\a_{\afrak} V  \\
L @>\xi>> \afrak\\
\end{CD}$$ and
$$\begin{CD}
L @>s>> E @>f>> \afrak\\
@| @. @| \\
L @>s'>>  E' @>f'>> \afrak.\\
\end{CD}$$We also have the analogous condition for M3a
\begin{description}
    \item[S1$\b$.] $\b_{\afrak}\circ\xi=\xi\circ\b_{L}$ and
    \item[S2$\b$.] $h\circ s=h'\circ s'$
\end{description}or, as above, in commutative diagram form
$$\begin{CD}
L @>\xi>> \afrak \\
@V\b_{L} VV @VV\b_{\afrak} V  \\
L @>\xi>> \afrak\\
\end{CD}$$ and
$$\begin{CD}
L @>s>> E @>h>> \afrak\\
@| @. @| \\
L @>s'>>  E' @>h'>> \afrak.\\
\end{CD}$$
The morphism $\f$ has to be subjected to these four conditions,
S1$\a,\b$ and S2$\a,\b$, for it to be a strong equivalence.
\subsection{Existence of extensions}
    So far we have shown how the 2-cocycle-like bilinear maps $g$ (that is, elements
    $g\in\Alt^2_\o(L,\afrak;\U)$ such that (\ref{eq:jacobi_gC}) holds) satisfying
    (\ref{eq:g'g}) corresponds in a \mbox{one-to-one}
    fashion to weak equivalence classes of extensions. We now show the existence of an extension $E$ of $L$ by $\afrak$
    given the following set of data.
    \begin{description}
        \item[Data A.] A function $f:L\oplus\afrak \to \afrak$ such that
        $f(0,a)=\a_{\afrak}(a)$.
        \item[Data B.] A function $h:L\oplus\afrak\to \afrak$ such that
        $h(0,a)=\b_{\afrak}(a)$.
        \item[Data C.] An alternating bilinear form $g$ satisfying (\ref{eq:jacobi_gC})
        and also
                    $$g(\a_L(x),\a_L(y))=h\big (\a_L(\langle
                    x,y\rangle_L),f(\langle x,y\rangle_L,g(x,y))\big ).$$
    \end{description}
    Now we put $E:=L\oplus \afrak$ and
    $$\langle (l,a),(l', a')\rangle_E\overset{\text{def}}{=}(\langle
    l,l' \rangle_L,g(l,l')).$$ In addition we choose
    the canonical section $s(x)=(x,0)$ and define
    $\pr$ and $\iota$ to be the natural projection and
    inclusion respectively:
    $$\pr:E\to L,\qquad\pr(x,a)=x;$$
    $$\iota:\afrak\to E,\qquad \iota(a)=(0,a).$$
    We also define $\o_E$ by
    \begin{align*}
    &\o_E(s(x)+\iota(a),s(y)+\iota(b))\circ s=s\circ \o_L(x,y),\qquad\text{
    and}\\
    &\o_E(s(x)+\iota(a),s(y)+\iota(b))\circ\iota=\iota\circ\o_\afrak(a,b)
    \end{align*}
    for $(x,y)\in D_{\o_L}$, $(a,b)\in D_\afrak$ and
    $(s(x)+\iota(a),s(y)+\iota(b))\in D_{\o_E}$. First note that the definition of the bracket can be written
    in the usual form $$\langle s(x),s(y)\rangle_E=s\langle x,y\rangle_L+\iota \circ g(x,y).$$
    This gives
    \begin{align*}
         \langle (l,a),(l',a')\rangle_E&=\big (\langle
         l,l'\rangle_L,g(l,l')\big )=\\
         &=s\langle l,l'\rangle_L+\iota\circ g(l,l')=\\
         &=s\circ
         \o_L(l,l')\langle
         l',l\rangle_L+\o_E(s(l),s(l'))\circ\iota\circ g(l',l)=\\
         &=\o_E(s(l),s(l'))\circ s\langle
         l',l\rangle_L+\o_E(s(l),s(l'))\circ\iota\circ g(l',l)=\\
         &=\o_E(s(l),s(l'))\circ \big (s\langle l',l\rangle_L+\iota\circ
         g(l',l)\big )
    \end{align*}which amounts to
    $$\langle (l,a),(l',a')\rangle_E=\o_E(s(l),s(l'))\big (\langle (l',a'),(l,a)\rangle_E\big ).$$
     Furthermore we define
    \begin{align*}
        &\a_E(x,a)\overset{\text{def}}{=}(\a_L(x),f(x,a))\\
        &\b_E(x,a)\overset{\text{def}}{=}(\b_L(x),h(x,a)).
    \end{align*}
    That $\a_E$ satisfies the $\b$-twisting condition follows from the following calculation.
    First,
    \begin{align*}
        \langle\a_E(x,a),\a_E(y,b)\rangle_E&=\langle
        (\a_L(x),f(x,a)),(\a_L(y),f(y,b))\rangle_E=\\
        &=\big (\langle\a_L(x),\a_L(y)\rangle_L,g(\a_L(x),\a_L(y))\big ).
    \end{align*}  Secondly,
    \begin{align*}
        &\b_E\circ\a_E\langle
        (x,a),(y,b)\rangle_E=\b_E\circ\a_E(\langle x,y\rangle_L,g(x,y))=\\
        &=\b_E(\a_L\langle x,y\rangle_L,f(\langle x,y\rangle_L,g(x,y)))=\\
        &=\Big(\b_L\circ\a_L\langle
        x,y\rangle_L,h\big(\a_L\langle x,y\rangle_L,f(\langle
        x,y\rangle_L,g(x,y))\big )\Big )
    \end{align*} and so comparing and using Data C gives us the
    result. Remember that $\b_E$ is only supposed to
    be a linear map by the definition of a qhl-algebra. The qhl-Jacobi
    identity is checked as
    \begin{align*}
        \big \langle \a_E(x,a),\langle (y,b),(z,c)\rangle_E\big
        \rangle_E&=\big \langle (\a_L(x),f(x,a)),(\langle
        y,z\rangle_L,g(y,z))\big \rangle_E=\\
        &=\big(\big \langle \a_L(x),\langle y,z\rangle_L\big
        \rangle_L,g(\a_L(x),\langle y,z\rangle_L)\big)
    \end{align*} and
    \begin{multline*}
        \b_E\big \langle (x,a),\langle (y,b),(z,c)\rangle_L\big
        \rangle_L=\b_E\big\langle (x,a),(\langle
        y,z\rangle_L,g(y,z))\big\rangle_E=\\
        =\b_E\big (\big \langle x,\langle y,z\rangle_L\big
        \rangle_L,g(x,\langle y,z\rangle_L)\big )=\\
        =\Big (\b_L\big \langle x,\langle
        y,z\rangle_L\big\rangle_L,h\big (\big \langle x,\langle
        y,z\rangle_L\big \rangle_L,g(x,\langle y,z\rangle_L)\big )\Big )
    \end{multline*}so we end up with
    \begin{align*}
        g(\a_L(x),\langle y,z\rangle_L)+h\big (\big \langle x,\langle
        y,z\rangle_L\big \rangle_L,g(x,\langle y,z\rangle_L)\big )
    \end{align*} since $L$ is a qhl-algebra. We know that $s$ is the canonical section $x\mapsto (x,0)$ and so by
    (\ref{eq:jacobi_gC}) the above expression is zero upon cyclic summation.

    That the diagram (\ref{diag:ses_ab1}) has exact rows is obvious from the definition of
    $\iota$ and $\pr$. Next we show that these maps are homomorphisms.
    \begin{align*}
       \pr\big(\langle
          (x,a),(y,b)\rangle_{E}\big)&=\pr\big(\langle
    x,y\rangle_L,g(x,y)\big)=\langle
    x,y\rangle_L=\\
     &=\langle\pr(x,a),\pr(y,b)\rangle_L,
    \end{align*}
    $$\langle \iota(a),\iota(b)\rangle_{E}=\langle
    (0,a),(0,b)\rangle_{E}=(0,0)=\iota(0)=\iota(\langle a,b\rangle_\afrak)$$
    since $\afrak$ was abelian. This shows that $\pr$ and $\iota$ are
    homomorphisms. In fact they are also qhl-algebra morphisms, because
    $$\pr\circ\a_E(x,a)=\pr(\a_L(x),f(x,a))=\a_L(x)=
\a_E\circ\pr(x,a)$$
    and
    $$\a_E\circ\iota(a)=\a_E(0,a)=(\a_L(0),f(0,a))=(0,\a_{\afrak}(a))=\iota\circ\a_{\afrak}(a)$$ and similary
    with $\b_E$. We have proved
    \begin{thm}\label{thm:homext2C} Suppose $(L,\Mor(L))$ and $(\afrak,\Mor(\afrak))$ are qhl-algebras with $\afrak$
    abel\-ian and put $E=L\oplus\afrak$.
    Then for every $g\in\Alt^2_{\o_L}(L,\afrak;\U)$ and every linear map
        $f:L\oplus\afrak\to\afrak$ and $h:L\oplus\afrak\to \afrak$ such that
        \begin{align}
            f(0,a)=\a_\afrak(a) \quad\text{for } a\in\afrak, \label{eq:cfvsaC}\\
            h(0,a)=\b_\afrak(a) \quad\text{for } a\in\afrak, \label{eq:cfvsb}
        \end{align}
        \begin{align}\label{eq:cbothgsandfsC}
            g(\a_L(x),\a_L(y))=h\big (\a_L(\langle
            x,y\rangle_L),f(\langle x,y\rangle_L,g(x,y))\big )
        \end{align}
    and
        \begin{multline}\label{eq:ctwococycle2C}
            \cs_{(x,a),(y,a),(z,c)}\,\o_E((z,c),(x,a))\circ\Big (\i\circ g(\a_L(x),\langle
            y,z\rangle_L)+\\
            +\i\circ h\big (\big \langle x,\langle y,z\rangle_L\big
            \rangle_L,g(x,\langle y,z\rangle_L)\big )\Big )=0,
        \end{multline} for $x,y,z\in L$ and $\big ((z,c),(x,a)\big ),\big ((x,a),(y,b)\big ),
        \big ((y,b),(z,c)\big )\in
        D_{\o_E}$, the linear direct sum $(E,\Mor(E))$ with morphisms $\alpha_E,\beta_E, \o_E$ is a qhl-algebra
    central extension of $(L,\Mor(L))$ by $(\afrak,\Mor(\afrak))$.
\end{thm}
\begin{remark}
Let $L$ be a qhl-algebra and $\afrak$ abelian.  We seek central
extensions of $L$ by
 $\afrak$. Put $E=L\oplus\afrak$ and let $$\langle \cdot,\cdot \rangle_E=\big (\langle\cdot,\cdot
 \rangle_L,g(\cdot,\cdot)\big)$$ for some bilinear $g:L\times L\to\afrak$. In addition to this pick the
 canonical section $s:L\to E, x\mapsto (x,0)$. Note that the definition of $\langle \cdot,\cdot\rangle_E$ is
 compatible with the choice of $s$. We want to endow $E$ with a qhl-algebra structure, that is, we
 want to find appropriate maps $\a_E$, $\b_E$ and $\o_E$ satisfying the necessary conditions. This means finding
 $f$ and $h$ such that $f(0,a)=\a_\afrak(a)$ and
 $h(0,a)=\b_\afrak(a)$. Define $\o_E$ by the relations
 $$\o_E(s(x)+\iota(a),s(y)+\iota(b))\circ s=s\circ \o_L(x,y)$$ and
    $$\o_E(s(x)+\iota(a),s(y)+\iota(b))\circ\iota=\iota\circ\o_\afrak(a,b).$$

 We let $\iota$ denote the canonical injection
$\iota:\afrak \hookrightarrow E, a\mapsto (0,a)$ and make the
 general ansatz
 \begin{align*}
    \iota\circ f(l,a)&=\big (0,\a_\afrak(a)+F(l)\big )\\
    \iota\circ h(l,a)&=\big (0,\b_\afrak(a)+H(l)\big ),
 \end{align*} for $F,H:L\to\afrak$ linear. A simple calculation shows that
 \begin{align*}
    \a_E(l,a)&=\big (\a_L(l),\a_\afrak(a)+F(l)\big )\\
    \b_E(l,a)&=\big (\b_L(l),\b_\afrak(a)+H(l)\big ).
 \end{align*}Furthermore, note that
 $$\big \langle s(x),s\langle y,z\rangle_E\big \rangle_E=\big (\big\langle x,\langle y,z\rangle_L\big
 \rangle_L,g(x,\langle y,z\rangle_L)\big),$$ and
 $$\iota\circ h\big \langle s(x),s\langle y,z\rangle_E\big \rangle_E=\big (0,\b_\afrak(g(x,\langle
 y,z\rangle_L))+H(\big \langle x,\langle y,z\rangle_L\big\rangle_L)\big ).$$ Putting this together to form the
 qhl-Jacobi identity we see
 \begin{align*}
    \cs_{(x,a),(y,a),(z,c)}\,\o_\afrak(c,a)\Big (&g\big (\a_L(x),\langle y,z\rangle_L\big)+\\
    &+\b_\afrak(g(x,\langle y,z\rangle_L))+
    H(\langle x,\langle y,z\rangle_L\rangle_L)\Big )=0,
 \end{align*} or, assuming $\b_\afrak=\id_\afrak$,
 \begin{align*}
    \cs_{(x,a),(y,a),(z,c)}\,\o_\afrak(c,a)\Big (g\big ((\a_L+\id_L)(x),\langle y,z\rangle_L\big)+
    H(\langle x,\langle y,z\rangle_L\rangle_L)\Big )=0.
 \end{align*}Notice the similarity with the corresponding hom-Lie algebra identity. In addition we must also
 have
 \begin{align*}
    \iota\circ g(\a_L(x),\a_L(y))&=\iota\circ h\big (\a_L\langle x,y\rangle_L,f\langle s(x),s(y)\rangle_E\big )=\\
    &=\big (0,\b_\afrak\circ f\langle s(x),s(y)\rangle_E+H(\a_L\langle x,y\rangle_L)\big )
 \end{align*} and so, if $\b_\afrak=\id_\afrak$,
  \begin{align*}
    g(\a_L(x),\a_L(y))&=f\langle s(x),s(y)\rangle_E+H(\a_L\langle x,y\rangle_L)=\\
    &=f(\langle
    x,y\rangle_L,g(x,y))+H(\a_L\langle x,y\rangle_L)=\\
    &=\a_\afrak(g(x,y))+F(\langle x,y\rangle_L)+H(\a_L\langle x,y\rangle_L).
 \end{align*}
\end{remark}
\begin{example}[Example \ref{ex:homLieExtC}, continued]
    Taking $L$ and $\afrak$ to be hom-Lie algebras, meaning that $h\circ\iota|_{\afrak}=\id_\afrak$ and
    $h\circ s|_{L}=0$, we see
    that we recover Theorem 7 from our previous paper
    {\rm\cite{HartLarsSilv1C}}. \qed
\end{example}
\begin{example}[Example \ref{ex:colorLieExt}, continued]
   Consider two color Lie algebras $L$ and $\afrak$ with the same
   grading group $\Gamma$ and the same commutation factor
   $\varepsilon$. This means that (see Example \ref{ex:colorLieExt})
   $$\afrak=\bigoplus_{\gamma\in\Gamma} \afrak_{\gamma}\qquad\text{and}\qquad
     L=\bigoplus_{\gamma\in\Gamma} L_{\gamma}.$$ Forming the vector
     space
   $$E=\bigoplus_{\gamma\in\Gamma} E_\gamma= \bigoplus_{\gamma\in\Gamma}
   (L_\gamma\oplus\afrak_\gamma)=L\oplus \afrak, $$we
   see immediately that
   $E$ is $\Gamma$-graded. We
   know from Theorem 9 and the deduction preceding it that we can
   endow this with a color structure as follows.

   From Examples \ref{ex:homLieExtC} and \ref{ex:colorLieExt} we see
   that $f\circ \iota|_\afrak=\id_\afrak$, $f\circ s|_L=0$
   and $h\circ\iota|_\afrak=\id_\afrak$, $h\circ s|_L=0$ and so {\rm(\ref{eq:cfvsaC})} and
   {\rm(\ref{eq:cfvsb})} are true. Taking $s: x\mapsto (x,0)$ and defining
   the product on $E$ by
    $$\langle (x,a),(y,b)\rangle_E=(\langle x,y\rangle_L,g(x,y))$$ for
    some $g\in\Alt^2_\varepsilon (L,\afrak;\U)$. That equation
    {\rm(\ref{eq:ctwococycle2C})}
    is satisfied we saw already in Example \ref{ex:colorLieExt}. Note that
    {\rm(\ref{eq:cbothgsandfsC})} becomes tautological. Hence we have a
    color central extension of $L$ by $\afrak$.

    Now $E$ is a color Lie algebra central extension of $L$ by
    $\afrak$. Note, however, that we have not constructed an explicit
    extension. What we have done is constructing an extension
    \emph{given} a $g\in\Alt^2_\varepsilon(L,\afrak;\U)$ satisfying
    {\rm(\ref{eq:jacobi_gC})} or
    its colored restriction. The existence of such a $g$ is \emph{not}
    guaranteed in general. See Scheunert {\rm\cite{Scheunert1C}} Proposition
    5.1 for a result that emphasizes this. In our setting this
    proposition implies that $H^2_{\varepsilon,-}(L,\afrak;\U)=\{0\}$
    and so there
    are no non-trivial central extensions. The actual construction of
    extensions, qualifying to finding 2-cocycles, is a highly
    non-trivial task.

    We now specialize the above slightly to one-dimensional central
    extensions with $\afrak=\k$. The $\k$-space $\k$ comes with a
    natural $\Gamma$-grading as
    $$ \k=\bigoplus_{\gamma\in\Gamma} K_\gamma, \qquad \text{where }
    K_0=\k,\quad K_\gamma=\{0\},\quad \text{ for } \gamma\neq 0.$$ Then we have
    a product on
    $$E=\bigoplus_{\gamma\in\Gamma}
    E_\gamma=\bigoplus_{\gamma\in\Gamma} (L_\gamma\oplus K_\gamma)$$
    defined by
    $$\langle (x,a),(y,b)\rangle_E=(\langle
    x,y\rangle_L,g(x,y)),$$where $g:L\times L\to \k$ is the
    $\k$-valued 2-cocycle.  \qed
\end{example}
\begin{remark}
    It would be of interest to develop a theory for quasi-hom-Lie
    algebra extensions of one qhl-algebra by another qhl-algebra,
    and apply it to get qhl-algebra extensions of the Virasoro algebra
     by a Heisenberg algebra {\rm\cite{JakobLee}}.
\end{remark}
\subsection{Central extensions of the $(\alpha,\beta,\o)$-deformed loop algebra}
Now form the vector space
$$\gc=\gh\oplus \k\cdot\cbf$$ for a ''central element'' $\cbf$ and pick the section
$s:\gh\to\gc, x\otimes t^n\mapsto (x\otimes t^n,0)$. Define a
bilinear product $\langle\cdot,\cdot\rangle_{\gc}$ on $\gc$ as
\begin{itemize}
    \item $\langle x\otimes t^n,\cbf\rangle_{\gc}=\langle \cbf,x\otimes t^n\rangle_{\gc}=0$, and
    \item $\langle x\otimes t^n,y\otimes t^m\rangle_{\gc}=\langle x,y\rangle_\g\otimes t^{n+m}+g(x\otimes t^n,y\otimes
    t^m)\cdot\cbf$, for a 2-cocycle-like bilinear map $g:\gh\times\gh\to\k$.
\end{itemize} It is easy to check that this definition is compatible with the section $s$.
Define, furthermore,
\begin{align*}
    &\a_{\gc}(x\otimes t^n+a\cdot\cbf)=\a_{\gh}(x\otimes t^n)+a\cdot\cbf\\
    &\b_{\gc}(x\otimes t^n+a\cdot\cbf)=\b_{\gh}(x\otimes t^n)+a\cdot\cbf,\quad\text{and}\\
    &\o_{\gc}(x\otimes t^n+a\cdot\cbf,y\otimes t^m+b\cdot\cbf)=\o_{\gh}(x\otimes t^n,y\otimes t^m)+\id.
\end{align*}We are now going to investigate when this 1-dimensional ''central extension'' $\gc$ of $\gh$
can be given the structure of a qhl-algebra.

The $\o$-skew symmetry of $\langle\cdot,\cdot\rangle_{\gc}$ and
$\b$-twisting of $\a_{\gc}$ is an easy check. To check the
$\o$-skew-symmetry condition observe that the product can be
written as, where we have put $u=x\otimes t^n$ and $v=y\otimes
t^m$ to simplify notation:
\begin{multline*}
    \langle u,v\rangle_{\gc}=\langle
    s(u),s(v)\rangle_{\gc}
    =s\langle u,v\rangle_{\gh}+\iota\circ
    g(u,v)=\\
    =s\circ\o_{\gh}(u,v)\langle v,u\rangle_{\gh}+\o_{\gc}(s(u),s(v))
    \circ\iota\circ g(v,u)=\\
    =\o_{\gc}(s(u),s(v))\circ s\langle
    v,u\rangle_{\gh}+\o_{\gc}(s(u),s(v))\circ\iota\circ g(v,u)=\\
    =\o_{\gc}(s(u),s(v))\circ\big (
    s\langle
    v,u\rangle_{\gh}+\iota\circ
    g(v,u)\big )=\o_{\gc}(s(u),s(v))\langle v,u\rangle_{\gc},
\end{multline*}and the $\beta$-twisting:
\begin{multline*}
    \langle
    \a_{\gc}(u+a\cdot\cbf),\a_{\gc}(v+b\cdot\cbf)\rangle_{\gc}=\langle
    \a_{\gh}(u)+a\cdot\cbf,\a_{\gh}(v)+b\cdot\cbf\rangle_{\gc}=\\
    =\langle
    \a_{\gh}(u),\a_{\gh}(v)\rangle_{\gh}=\b_{\gh}\circ\a_{\gh}\langle
    u,v\rangle_{\gh}=\b_{\gc}\circ\a_{\gc}\langle
    u+a\cdot\cbf,v+b\cdot\cbf\rangle_{\gh}.
\end{multline*}
     For the qhl-Jacobi identity we have first
\begin{multline}
    \big \langle \a_{\gc}(x\otimes t^n),\langle y\otimes t^m,z\otimes t^l\rangle_{\gc}\big
    \rangle_{\gc}=\\
    =\big \langle \a_{\g}(x)\otimes t^n,\langle y,z\rangle_\g\otimes t^{m+l}+g(y\otimes t^m,z\otimes
    t^l)\cdot\cbf\big \rangle_{\gc}=\nonumber\\
    =\big \langle \a_{\g}(x)\otimes t^n,\langle y,z\rangle_\g\otimes t^{m+l}\big \rangle_{\gc}=\nonumber\\
    =\big \langle
    \a_\g(x),\langle y,z\rangle_\g\big \rangle_\g\otimes t^{n+m+l}+g(\a_\g(x)\otimes t^n,\langle y,z\rangle_\g\otimes t^{m+l})\cdot\cbf.
\end{multline}Secondly,
\begin{multline}
    \b_{\gc}\big\langle x\otimes t^n,\langle y\otimes t^m,z\otimes t^l\rangle_{\gc}\big \rangle_{\gc}=\notag\\
    =\b_{\gc}\big (\big \langle x,\langle y,z\rangle_\g\big \rangle_\g\otimes t^{n+m+l}+
    g(x\otimes t^n,\langle
    y,z\rangle_\g\otimes t^{m+l})\cdot\cbf\big )=\notag\\
    =\b_\g\big (\big\langle x,\langle y,z\rangle_\g\big \rangle_\g\big )\otimes t^{n+m+l}+
    g(x\otimes t^n,\langle
    y,z\rangle_\g\otimes t^{m+l})\cdot\cbf
\end{multline} and so combining and summing up cyclically, using that $\g$ is a qhl-algebra, we end up with
\begin{align}\label{eq:2cocycLoop}
    \cs_{(x,n),(y,m),(z,l)}\, g((\a_\g+\id_\g)(x)\otimes t^n,\langle y,z\rangle_\g\otimes
    t^{m+l})=0.
\end{align}
Now to do this a little more explicit and more in tune with the
classical Lie algebra case (\cite{FLMC}) we construct the product
on $\gc$ a bit differently. We assume
$\o_{\g}=\o_{\gh}=\o_{\gc}=-1$, that is, that the product is
skew-symmetric, and take a homogeneous $\sigma$-derivation $D$ on
$\k[t,t^ {-1}]$, homogeneous in the sense that $\sigma$ is a
homogeneous mapping $t\mapsto qt$, for $q\in\k^*$, the
multiplicative group of non-zero elements of $\k$.\footnote{One
would be tempted to try a more general
  $\sigma$-derivation with $\sigma(t)=qt^s$ as in section \ref{sec:nonlinC}, but
  the following construction seems only to work with $s=1$.} Explicitly we have
  $$D=\eta t^{-k}\frac{\id-\sigma}{1-q}$$ and leading to
$$ D(t^n)=\eta\{n\}_qt^{n-k}.$$
Take a bilinear form $B(\cdot,\cdot)$ on $\g$ and factor the
2-cocycle-like bilinear map $g$ as
$$g(x\otimes t^n,y\otimes t^m)=B(x,y)\cdot (D(t^n)\cdot t^m)_0,$$
where the notation $(f)_0$ is the zeroth term in the Laurent
polynomial $f$ or, put differently, $t$ times the residue
$\Res(f)$. The above trick to factor the 2-cocycle (in the Lie
algebra case) as $B$ times a ''residue'' is apparently due to Kac
and Moody from their seminal papers where they introduced what is
now known as Kac-Moody algebras, \cite{KacC} and \cite{MoodyC},
respectively. This means that
$$(D(t^n)\cdot t^m)_0=\eta\{n\}_q\delta_{n+m,k}.$$
Calculating the 2-cocycle-like condition (\ref{eq:2cocycLoop}) now
leads to
$$\cs_{(x,n),(y,m),(z,l)}\,(\eta\cdot\{n\}_q\cdot\delta_{n+m+l,k})\cdot
B((\alpha_\g+\id)(x),\langle y,z\rangle_\g)=0,$$and for
$\alpha_\g=\id,\eta=1,q=1$ and $k=0$ we retrieve the classical
2-cocycle discovered by Kac and Moody. Notice, however, that in
the Lie algebra case it is assumed that $B$ is symmetric and
$\g$-invariant, this leading to a nice 2-cocycle identity unlike
the one we have here. What we thus obtained by the preceding
factorization is a \mbox{$(\alpha,\beta,-1)_q$-deformed,}
one-dimensional central extension of the (Lie) loop algebra, where
the $q$-subscript is meant to indicate that we have $q$-deformed
the derivation on the Laurent polynomial as well as the underlying
algebra. Note that the 2-cocycle-like condition only depends on
the ''base algebra'' $\g$.

\section*{Acknowledgements}
We would like to express our gratitude to Jonas Hartwig for
valuable comments and insights and to Hans Plesner Jakobsen for
bringing to our attention the reference \cite{JakobLee}. The first
author stayed at the Mittag-Leffler Institute, Stockholm, during
the last phase of writing this paper in January-March 2004; the
second author stayed there in September-October 2003 and May-June
2004. Very warm thanks go out to Mittag-Leffler Institute for
support and to the staff and colleagues present for making it a
delightful and educating stay.

Some results appearing in this paper were reported by the second
author at the Non-commutative Geometry Workshop I, Mittag-Leffler
Institute, Stockholm, September, 2003, and by the first author at
the 3'rd {\"O}resund Symposium in Non-commutative Geometry and
Non-commutative Analysis, Lund January, 2004.


\begin{thebibliography}{99}
\bibitem{AizSato1C} Aizawa, N., Sato, H.-T., \emph{$q$-deformation of the
Virasoro algebra with central extension}, Phys. Lett. B, 256
(1999), no. 2, 185--190.
%
\bibitem{BlochS1C} Bloch S., \emph{Zeta values and differential operators on the
circle}, J. Algebra, 182 (1996), 476--500.
\bibitem{ChaiIsLukPopPresC} Chaichian, M., Isaev, A. P.,
Lukierski, J., Popowicz, Z., \mbox{Pre\v snajder, P.,}
\emph{\mbox{$q$-deformations} of Virasoro algebra and conformal
dimensions}, Phys. Lett. B 262 (1991), no. 1, 32-38.
\bibitem{ChaiKuLukC} Chaichian, M., Kulish, P.,
Lukierski, J., \emph{$q$-deformed Jacobi identity,
\mbox{$q$-oscillators} and $q$-deformed infinite-dimensional
algebras}, Phys. Lett. B 237 (1990), no. 3-4, 401--406.
\bibitem{ChaiPres3C} Chaichian, M., Pre\v snajder, P.,
\emph{$q$-Virasoro algebra, $q$-conformal dimensions and free
$q$-superstring}, Nuclear Phys. B 482 (1996), no. 1-2, 466--478.
\bibitem{ChungWSC} Chung, W.-S., \emph{Two parameter deformation of Virasoro
algebra}, J. Math. Phys. 35 (1994), no. 5, 2490--2496.
\bibitem{CurtrZachos1C} Curtright, T. L., Zachos, C. K.,
\emph{Deforming maps for quantum algebras}, Phys. Lett. B 243
(1990), no. 3, 237--244.
\bibitem{QuantFieldStrings} Deligne, P., Etingof, P., Freed, D.S.,
Jeffrey, L.C., Kazhdan, D., \mbox{Morgan, J.W.,} Morrison, D.R.,
Witten, E., (Eds) \emph{Quantum Fields and Strings: A Course for
Mathematicians}, 2 vol., Amer. Math. Soc., 1999.\\
(ISBN: 0-8218-2014-1).
\bibitem{DiFranMathiSenC} Di
Francesco, P., Mathieu, P., S\'en\'echal, D., \emph{Conformal
Field Theory}, Springer Verlag,
1997, 890 pp.\\
(ISBN: 0-387-94785-X).
 \bibitem{FLMC} Frenkel, I., Lepowsky, J., Meurman, A., \emph{Vertex Operator Algebras and the Monster},
    Academic Press, 1988, 508 pp.\\
    (ISBN: 0-12-267065-5).
\bibitem{Fuchs1C} Fuchs, J., \emph{Affine Lie Algebras and Quantum
Groups}, Cambridge University Press, 1992, 433 pp.\\
(ISBN: 0 521 41493 4).
\bibitem{Fuchs2C} Fuchs, J., \emph{Lectures on Conformal Field
Theory and Kac-Moody Algebras}, \texttt{hep-th/9702194}, February
1997.
\bibitem{FuksC} Fuks, D.B., \emph{Cohomology of Infinite-Dimensional Lie
Algebras}, Plenum Publishing Corp., 1986\\
(ISBN: 0-306-10990-5).
\bibitem{HartLarsSilv1C} Hartwig, J.T., Larsson D., Silvestrov S.D.,
  \emph{Deformations of Lie algebras using $\sigma$-derivations}, Preprints
  in Mathematical Sciences 2003:32, LUTFMA-5036-2003, Centre for
  Mathematical Sciences, Department of Mathematics, Lund Institute of
  Technology, Lund University, 2003.
\bibitem{HelSil-bookC}
{\rm Hellstr{\"o}m L., Silvestrov, S.D.}, \emph{Commuting Elements
in $q$-Deformed Heisenberg Algebras},
World Scientific, 2000, 256 pp.\\
 (ISBN: 981-02-4403-7).
\bibitem{JakobLee} Jakobsen, H.P., Lee, H.C.-W., \emph{Matrix
Chain Models and Kac-Moody Algebras}, To appear in "Vertex
Operators and Kac-Moody algebras" Proceedings Madras 2002,
Contemp. Math., Amer. Math. Soc.
\bibitem{KacC} Kac, V.G.,
\emph{Simple irreducible graded Lie algebras of finite growth},
Math. USSR Izv. 2, 1968, 1271--1311.
\bibitem{KacRadulC} Kac, V.G., Radul, A., \emph{Quasifinite highest weight modules over
the Lie algebra of differential operators on the circle}, Commun.
Math. Phys. 157, 1993, 429--457.
\bibitem{KacRainaC} Kac, V.G., Raina, A.K., \emph{Highest weight representations of infinite-dimensional Lie
algebras}, World Scientific, 1987, 145 pp. \\
(ISBN: 9971-50-395-6).
\bibitem{Kassel1C} Kassel, C., \emph{Cyclic homology of differential operators,
the Virasoro algebra and a $q$-analogue}, Commun. Math. Phys. 146
(1992), 343-351.
\bibitem{KhesinLyubRogerC} Khesin, B., Lyubashenko, V., Roger, C.,
\emph{Extensions and contractions of Lie algebra of
$q$-Pseudodifferential symbols on the circle}, J. Func. Anal. 143
(1997), 55-97.
\bibitem{LiC} Li, W.-L., \emph{$2$-Cocycles on the algebra of differential operators},
J. Alg. 122 (1989), 64-80.
\bibitem{MoodyC} Moody, R.V., \emph{A new class of Lie algebras},
J. Algebra 10, 1968, 211--230.
\bibitem{Passman1C} Passman, D.S., \emph{Simple Lie color algebras of Witt type}, J. Algebra 208, (1998),
698--721.
\bibitem{PolychronakosC} Polychronakos, A. P., \emph{Consistency conditions
and representations of a $q$-deformed Virasoro algebra}, Phys.
Lett. B 256 (1991), no. 1, 35--40.
\bibitem{SatoHT1C} Sato, H.-T., \emph{Realizations of $q$-deformed Virasoro algebra},
Progress of Theoretical Physics, 89 (1993), no. 2, 531-544.
\bibitem{SatoHT2C} Sato, H.-T., \emph{$q$-Virasoro operators from an
analogue of the Noether currents}, Z. Phys. C 70 (1996), no. 2,
349--355.
\bibitem{Scheunert1C} Scheunert, M., \emph{Introduction to the cohomology of Lie superalgebras and some
    applications}, Research and Exposition in Mathematics, Volume 25 (2002), 77--107.
\bibitem{ScheunertZhangC} Scheunert, M., Zhang, R.B.,
\emph{Cohomology of Lie superalgebras and their generalizations},
J. Math. Phys. 39, 1998, 5024--5061.
\bibitem{Su1C} Su, Y., \emph{$2$-Cocycles on the Lie algebras of generalized
differential operators}, Comm. Alg. 30 (2) (2002), 763-782.
\end{thebibliography}
\end{document}